\documentclass[11pt]{article}

\usepackage[T1,OT1]{fontenc}
\usepackage[margin=1in]{geometry}
\usepackage{amsmath,amssymb,amsthm,mathtools}
\usepackage{booktabs}
\usepackage{enumitem}
\setlist[itemize]{label=$\bullet$}
\usepackage{microtype}
\usepackage{xcolor}
\usepackage{graphicx}
\usepackage[colorlinks=true,linkcolor=blue!50!black,citecolor=blue!50!black,urlcolor=blue!50!black]{hyperref}

\newtheorem{theorem}{Theorem}[section]
\newtheorem{lemma}[theorem]{Lemma}
\newtheorem{proposition}[theorem]{Proposition}
\newtheorem{corollary}[theorem]{Corollary}
\newtheorem{heuristic}[theorem]{Heuristic}
\theoremstyle{definition}
\newtheorem{definition}[theorem]{Definition}

\theoremstyle{remark}
\newtheorem{remark}[theorem]{Remark}

\newcommand{\Zp}{\mathbb{Z}_p}
\newcommand{\Qp}{\mathbb{Q}_p}
\newcommand{\N}{\mathbb{N}}
\newcommand{\abs}[1]{\lvert #1 \rvert}
\newcommand{\norm}[1]{\lVert #1 \rVert}
\DeclareMathOperator{\vp}{v}
\newcommand{\B}{\mathrm{B}}
\newcommand{\doi}[1]{\href{https://doi.org/#1}{doi:#1}}
\newcommand{\ip}[2]{\langle #1, #2 \rangle}

\title{Certified $p$-adic normal forms:\\ an Eliasson--Vey theorem and zero sets of orbits}
\author{Kirill Terskii\\ \texttt{t02260049@gse.cs.msu.ru}}
\date{}

\begin{document}

\maketitle

\begin{abstract}
I use the Smith form of the homological operator to construct partial Poincar\'e--Dulac normal forms near an equilibrium of a vector field on $\mathbb Z_p^d$ with coefficients divisible by $p$ (by $4$ if $p=2$). Retaining the directions whose removal would cost more than a fixed number of digits gives convergence on an explicit ball for any linear part. For polynomial fields, an exact integer check of the conjugacy defect certifies the computed flow to $K$ digits in rescaled coordinates. When the quadratic parts of an integrable Hamiltonian system at a rank-$0$ point span a Cartan subalgebra, the homological divisors in symplectic eigencoordinates are integers. This lets me use Newton iteration to obtain an exact $p$-adic Eliasson--Vey normal form at every prime and for every Williamson type. In exact linear normal coordinates, a polynomial restricted to an orbit gives a sum of exponentials, to which van der Poorten's bound of $2s-3$ zeros applies for a nonzero sum with $s\ge2$ distinct exponents. For integer orbits of polynomial maps, I combine exact modular orbits, Strassmann bounds and archimedean estimates in an adelic sieve. At every prime, increasing the modular depth eventually excludes a time class exactly when the interpolation along its center-manifold shadow has no zero in $\mathbb Z_p$.
\end{abstract}

\section{Introduction}

The Poincar\'e--Dulac construction simplifies a vector field near an equilibrium by solving a homological equation at each degree \cite{arnold}. For a diagonal linear part this means dividing by numbers of the form $\ip m\lambda-\lambda_i$, losing as many $p$-adic digits as their valuations. In general the eigenvalues can understate the loss. A $k\times k$ Jordan block $3I+J$ over $\mathbb Z_3$, for example, has only the eigenvalue $3$, yet its Smith exponents are $0,\ldots,0,k$, so inverting it costs $k$ digits.

I use the Smith form of the homological operator to decide which directions to remove, retaining those whose removal costs more than $\sigma$ digits. Rescaling about the equilibrium supplies the digits needed at each degree. For a field of valuation $c$, Theorem~\ref{thm:threshold-convergence} gives convergence on the ball of radius $p^{-r}$ for $r=1+\max(\sigma-c,0)$, including resonant and non-semisimple linear parts, and bounds the tails of both series. At finite precision I evaluate the conjugacy defect directly on the computed polynomials. Bounding its valuation gives agreement of the flows modulo $p^{K+r}$ on the ball, or $K$ digits in the rescaled coordinates (Theorem~\ref{thm:certificate}). If the field is a power series, I also need a lower bound on the valuations of the omitted coefficients.

Herman and Yoccoz proved a non-Archimedean analog of Siegel's linearization theorem under a Diophantine condition \cite{hermanyoccoz}, and Lindahl determines the size of the linearization disk for quadratic maps in one variable \cite{lindahl}. Jenkins and Spallone treat formal and analytic flows and the classification of maps tangent to the identity \cite{jenkinsspallone}. The threshold has a parallel in linear differential equations: while the $p$-adic Fuchs theorem uses a non-Liouville condition \cite{kedlaya}, Kedlaya's Liouville partitions group exponents with integer or Liouville differences. For a differential module on a disk satisfying the Robba condition, those blocks give a direct sum decomposition \cite[Definition~3.4.4 and Remark~3.4.15]{kedlaya15}.

The Smith exponents are the Hodge slopes of the homological operator, and its eigenvalue valuations are the Newton slopes. The classical comparison \cite{kedlaya,elsheikhgiesbrecht} bounds the cost of an individual division here. If the linear part diagonalizes in a basis $S$, Proposition~\ref{prop:smith-semisimple} gives $\abs{s_j-v_j}\le(n+1)\kappa(S)$ for the ordered slopes in degree $n$, with $\kappa(S)$ measuring the failure of the eigenbasis to be integral. Algebraic eigenvalues with an integral eigenbasis give an $O(\log n)$ bound, while the Liouville example shows how fast the exponents can grow without such hypotheses. Section~\ref{sec:formal-eigen} also treats logarithms of algebraic points in the kernel of reduction. At an odd prime and with an integral eigenbasis, exact group arithmetic gives the exponents and heights bound them linearly in the degree for $\mathbb G_m$ or quadratically for an elliptic curve. For $\mathbb G_m$, Yu's theorem improves the growth to a polynomial in $\log n$ with larger constants \cite{yu}.

If the quadratic parts of commuting Hamiltonians at a rank-$0$ point span a Cartan subalgebra, the divisors in symplectic eigencoordinates are integers with valuation at most $\log_p$ of the degree. This gives a Newton iteration for $\mathcal H\circ\Phi=\mathcal G(q)$ on an explicit ball (Theorem~\ref{thm:ev}). Crespo and Pelayo ask which Williamson types admit such a local form \cite[Question~5.14]{crespopelayoww}, and the rank-$0$ theorem covers every type and every prime under the Cartan condition in their later formulation \cite[Definition~3.2]{crespopelayocontrast}.

In linear normal coordinates, a polynomial restricted to an orbit becomes a sum of exponentials whose exponents come from the eigenvalues. For a nonzero sum with $s\ge2$ distinct exponents, van der Poorten and Schlickewei give a bound of $2s-3$ zeros that is uniform in the starting point \cite{vdpoorten,vdpoortenschlickewei}. In Section~\ref{sec:zeros} I ask for the exact integer times at which an orbit of a polynomial map meets a hypersurface \cite{ostafesha}. I start with the orbit modulo $p^k$. When a time class survives at this depth, the center-manifold shadow tells whether raising $k$ can ever exclude it. On a nilpotent cycle the functions on the shadow are just constants, and their valuations give the depth bound in Proposition~\ref{prop:nilpotent}. To estimate how often the orbit has an \'etale cycle modulo $p$, I use the random-map model of \cite{bghkst}.

Write $\vp$ for the $p$-adic valuation, taken as the minimum over entries of a vector or coefficients of a polynomial. Thus $\norm{x}=p^{-\vp(x)}$, and $\B(x,R)$ denotes the closed ball of radius $R$ centered at $x$. For a vector field $f$, set $c=\vp(f)$ and assume
\begin{equation}\label{eq:hyp}
c > \tfrac{1}{p-1}, \quad\text{that is,}\quad c \ge 1 \text{ for } p \ge 3, \qquad c \ge 2 \text{ for } p = 2.
\end{equation}
Time rescaling gives this condition for any field with $p$-integral coefficients after multiplication by $p$, or by $4$ when $p=2$. It is sharp for linear fields because the flow $e^{\lambda t}x$ of $f(x)=\lambda x$ is defined for all $t\in\Zp$ exactly when $\vp(\lambda)>1/(p-1)$. I use Legendre's formula $\vp(n!)=(n-s_p(n))/(p-1)\le(n-1)/(p-1)$, where $s_p(n)$ is the sum of the base-$p$ digits of $n$.

\section{Flows of \texorpdfstring{$p$}{p}-adic vector fields}\label{sec:flow}

Let $D = \sum_j f_j \partial_j$ be the derivation attached to $f$, $E_n = D^n x/n!$, applied componentwise to the coordinate functions, and $\varphi_t(x) = \sum_{n \ge 0} t^n E_n(x)$.

\begin{lemma}[{\cite[Lemma~2.1]{terskiiA}}]\label{lem:coeff}
$\vp(D^n g) \ge \vp(g) + nc$ for every $g \in \Zp[x]$, and $\vp(E_n) \ge nc - \vp(n!) \ge n\bigl(c - \tfrac{1}{p-1}\bigr) + \tfrac{1}{p-1}$ for $n\ge1$.
\end{lemma}

So the Lie series converges on $\Zp \times \Zp^d$, and since $D$ is a derivation it defines a flow: $\varphi_{s+t} = \varphi_t \circ \varphi_s$ and $\partial_t\varphi_t(x) = f(\varphi_t(x))$. Among the solutions given by power series in $t$ converging on $\Zp$, the map $t \mapsto \varphi_t(x)$ is the only one with $\varphi_0(x) = x$, because the coefficients of such a series are determined recursively. Everything holds verbatim for $f$ in the Tate algebra $\Zp\langle x\rangle^d$ of power series with coefficients tending to zero. By Lemma~\ref{lem:coeff} the coefficients of $E_n$, $n \ge 1$, are divisible by $p$, so each $\varphi_t$ is an isometry of $\Zp^d$, the identity plus a map that is $p^{-1}$-Lipschitz \cite[Section~2]{terskiiA}.

\begin{lemma}[{\cite[Theorem~3.2]{terskiiA}}]\label{lem:orbit-metric}
If $f(x) \ne 0$, then $\norm{\varphi_t(x) - \varphi_s(x)} = \norm{f(x)}\,\abs{t - s}$ for all $s, t \in \Zp$.
\end{lemma}

\section{Normal forms in any dimension}\label{sec:normal}

An approximate equilibrium can be certified before constructing the normal form. For $f\in\Zp[x_1,\ldots,x_d]^d$ and $a\in\Zp^d$ with $\det Df(a)\ne0$, put $m=\vp(f(a))$ and $e=\vp(\det Df(a))$. If $m>2e$, Hensel's lemma gives a simple zero in $\B(a,p^{-(m-e)})$, unique in the larger ball $\B(a,p^{-e-1})$. Indeed, the adjugate formula gives $\norm{Df(a)^{-1}}\le p^e$. On the larger ball the map $x\mapsto x-Df(a)^{-1}f(x)$ has Lipschitz constant at most $p^{-1}$ by the Taylor remainder estimate about $a$, and it moves $a$ by at most $p^{-(m-e)}$. It therefore has a fixed point in the smaller ball, and $Df$ stays invertible on the larger one.

Around a simple equilibrium $x^*$ in dimension one, with $\lambda = f'(x^*)$ and $f(x^* + y) = \lambda y + \sum_{k \ge 2} a_k y^k$, the Poincar\'e--Siegel coordinate $h(y) = y + \sum_{n \ge 2} h_n y^n$ with $h \circ \varphi_t = e^{\lambda t} h$ satisfies
\[
(n-1)\lambda\, h_n = -\sum_{k=2}^{n} (n-k+1)\, a_k\, h_{n-k+1}.
\]
Since $x^*\in\Zp$, all $a_k$ and $\lambda$ have valuation at least $c$. With $\ell = \vp(\lambda)\ge c$ the recursion gives $\vp(h_n) \ge (n-1)(c - \ell) - \vp((n-1)!)$, so $h$ converges on $\abs{y} < p^{-\rho}$ with $\rho = \ell - c + 1/(p-1)$, and $\vp(h_ny^n)\ge(n-1)(\vp(y)-\rho)+\vp(y)$. The coordinate gives an independent way to compute the flow and orbit times near $x^*$, through $e^{\lambda t}$ and the $p$-adic logarithm. In one dimension the divisors are $(n-1)\lambda$, so there is no small-divisor obstruction.

Let $x^*$ be an equilibrium and $A = Df(x^*)$. Fix $r \ge 1$, put $x = x^* + p^r z$, and
\[
\tilde f(z) = p^{-r} f(x^* + p^r z) = Az + \sum_{n \ge 2} \tilde f_n(z), \qquad \vp(\tilde f_n) \ge c + r(n - 1).
\]
As in the classical theory \cite{arnold}, I look for a change of coordinates $\Phi(z) = z + \sum_{n \ge 2} \Phi_n(z)$ and a field $g_{\mathrm{nf}}(z) = Az + \sum_{n \ge 2} R_n(z)$, with $\Phi_n, R_n$ homogeneous of degree $n$, such that $\Phi$ conjugates $g_{\mathrm{nf}}$ to $\tilde f$:
\begin{equation}\label{eq:conj}
D\Phi(z)\, g_{\mathrm{nf}}(z) = \tilde f(\Phi(z)).
\end{equation}
Comparing the terms of degree $n$ gives the homological equation
\begin{equation}\label{eq:homological}
L_A \Phi_n + R_n = G_n, \qquad L_A \Phi = D\Phi \cdot Az - A\Phi,
\end{equation}
\[
G_n = \Bigl[ \tilde f_{\ge 2}\bigl(z + \textstyle\sum_{k < n} \Phi_k\bigr) \Bigr]_n - \sum_{\substack{j + k = n + 1 \\ 2 \le j, k < n}} D\Phi_j \cdot R_k,
\]
where $G_n$ depends only on the terms of lower degree. If $R = 0$, then $\varphi_t(x) = x^* + p^r \Phi(e^{At} \Phi^{-1}((x - x^*)/p^r))$, which is the linearization. In general $\Phi$ conjugates the flow of $\tilde f$ to the flow $\psi_t$ of $g_{\mathrm{nf}}$.

Let $V_n$ be the free $\Zp$-module of homogeneous vector polynomials of degree $n$ in $d$ variables, with the basis of monomials $z^e \varepsilon_i$, where $\varepsilon_1, \dots, \varepsilon_d$ is the standard basis of $\Zp^d$. In this basis $L_A$ is a square matrix with entries in $\Zp$ (entries of $A$ times exponents). Its Smith normal form over $\Zp$ is \cite{newman}
\[
P L_A Q = \operatorname{diag}\bigl(p^{s_1}, \dots, p^{s_m}\bigr), \qquad P, Q \in \mathrm{GL}_m(\Zp),
\]
where $p^{\infty}=0$, and in practice the decomposition is computed modulo a working power of $p$ by Gaussian elimination with pivots of least valuation. Put $\gamma = P G_n$. Fix an integer threshold $\sigma\ge0$. Set $\xi_i=\gamma_i/p^{s_i}$ if $s_i\le\sigma$. For $s_i>\sigma$, set $\xi_i=0$ and keep $\gamma_i$ in the normal form. The divisions take place in $\Qp$. Theorem~\ref{thm:threshold-convergence} gives a sufficient radius for the quotients to be integral. Then
\[
\Phi_n = Q\xi, \qquad R_n = P^{-1} \textstyle\sum_{s_i > \sigma} \gamma_i \varepsilon_i
\]
solve \eqref{eq:homological}, since $P(L_A \Phi_n + R_n) = \operatorname{diag}(p^{s_i}) \xi + \sum_{s_i > \sigma} \gamma_i \varepsilon_i = \gamma$.

\begin{lemma}\label{lem:keep}
For every choice of the kept directions, the truncations of $\Phi$ and $g_{\mathrm{nf}}$ at degree $n$ satisfy \eqref{eq:conj} up to terms of degree greater than $n$. In particular, keeping a removable term in $R$ never makes the conjugacy inexact.
\end{lemma}

\begin{proof}
The terms of degree $k \le n$ of \eqref{eq:conj} are exactly the equations \eqref{eq:homological} for $k \le n$, which hold by construction.
\end{proof}

The threshold does not affect the formal identity. It trades the size of the domain, which shrinks with the divisors that are divided out, against the number of terms left in $g_{\mathrm{nf}}$. Exact zero Smith factors are always kept. At working precision $W>\sigma$, a zero factor modulo $p^W$ means $s_i\ge W$ and is kept. Since removable directions can also be kept, this is a partial normalization of Poincar\'e--Dulac type, not necessarily a normal form containing only resonant terms.

\begin{proposition}\label{prop:divisors}
Suppose that $A$ is diagonalizable over a finite extension of $\Qp$, with eigenvalues $\lambda_1, \dots, \lambda_d$. Then $L_A$ is diagonalizable with the eigenvalues $\ip{m}{\lambda} - \lambda_i$, $\abs{m} = n$. The number of directions with $s_i = \infty$ is the number of resonances $\ip{m}{\lambda} = \lambda_i$, and if there is none, $\sum_i s_i = \sum_{\abs{m} = n,\, i} \vp\bigl(\ip{m}{\lambda} - \lambda_i\bigr)$.
\end{proposition}

\begin{proof}
In an eigenbasis of $A$, $L_A(z^m \varepsilon_i) = (\ip{m}{\lambda} - \lambda_i)\, z^m \varepsilon_i$. The kernel of a diagonalizable operator is the eigenspace of $0$, and its dimension is the rank defect, which the Smith form counts. The determinant of $L_A$ is the product of its eigenvalues, and $\det(P L_A Q)$ equals it up to a unit, so the valuations add up.
\end{proof}

The elementary divisors thus carry the information of the classical small divisors without any eigenvalue being computed. This is useful because the eigenvalues usually lie in extensions of $\Qp$: complex pairs of the real theory have no special status here, and an irreducible characteristic polynomial is handled like any other. For a Jordan block the $s_i$ are not individual divisors: along a chain they add up. For example, the $k \times k$ matrix $3I + J$, with $J$ the nilpotent Jordan block (ones above the diagonal), has Smith exponents $0, \dots, 0, k$ over $\mathbb{Z}_3$, that is, elementary divisors $1, \dots, 1, 3^k$, while its only eigenvalue $3$ has valuation $1$. Zero eigenvalues produce infinitely many resonances: $\ip{m}{\lambda} = \lambda_i$ for every $m = m' + e_i$, where $e_i$ is the $i$-th unit exponent and $m'$ is supported on the zero directions, and the normal form keeps the corresponding series.

\subsection{Bounds on the Smith exponents}\label{sec:smith-bounds}

Write $L_n=L_A|_{V_n}$, $M_n=d\binom{n+d-1}{d-1}$ and $r_n=\operatorname{rank}_{\Qp}L_n$. Order its exponents as
\[
s_1(n)\le\cdots\le s_{r_n}(n)<\infty,
\qquad s_j(n)=\infty\quad(j>r_n).
\]
When $r_n>0$, put $\mu_n=s_{r_n}(n)$. The relevant bound is an upper bound on these valuations: division by $p^{s_i(n)}$ loses $s_i(n)$ digits. For integer matrices, when Smith exponents and eigenvalue valuations coincide is studied in \cite{elsheikhgiesbrecht}. Here the matrix is the homogeneous representation $L_n$ of $A$. In the language of \cite{kedlaya}, the $s_j(n)$ are the slopes of the Hodge polygon of $L_n$, and the valuations $v_1(n) \le \cdots \le v_{M_n}(n)$ of its eigenvalues $\ip{m}{\lambda} - \lambda_i$, the small divisors, with $\infty$ for zero, are the slopes of its Newton polygon. The Newton polygon lies above the Hodge polygon, with the same endpoint when $L_n$ is invertible:
\begin{equation}\label{eq:newton-hodge}
\sum_{j \le i} s_j(n) \le \sum_{j \le i} v_j(n) \quad (1 \le i \le M_n), \qquad \sum_{j \le M_n} s_j(n) = \sum_{j \le M_n} v_j(n) \ \text{ if } r_n = M_n.
\end{equation}
In particular $s_1(n) \le v_1(n)$: the smallest exponent never exceeds the valuation of the smallest divisor. For singular $L_n$, \eqref{eq:newton-hodge} follows from the invertible case applied to $L_n+p^u I$, letting $u\to\infty$.

\begin{proposition}[Minors]\label{prop:smith-minors}
Suppose $r_n>0$ and put $a=\min_{i,j}\vp(A_{ij})$. If $\Delta_j(n)$ is the minimum valuation of the non-zero $j\times j$ minors of $L_n$, with $\Delta_0(n)=0$, then
\begin{equation}\label{eq:smith-minors}
s_j(n)=\Delta_j(n)-\Delta_{j-1}(n)\ \ (1\le j\le r_n),\qquad
a\le s_j(n),\qquad
\mu_n\le\Delta_{r_n}(n)-(r_n-1)a.
\end{equation}
\end{proposition}

\begin{proof}
This is the description of elementary divisors by minors \cite{newman,kedlaya}: the ideal generated by the $j\times j$ minors is unchanged by multiplication on either side by an invertible integral matrix. For a Smith diagonal matrix its valuation is $s_1+\cdots+s_j$, which proves the first identity. Every entry of $L_n$ is divisible by $p^a$, so every finite Smith exponent is at least $a$, giving the last inequality in \eqref{eq:smith-minors}.
\end{proof}

When $L_n$ is invertible, \eqref{eq:newton-hodge} computes $\Delta_{M_n}(n) = \vp(\det L_n)$ exactly, and the last bound of \eqref{eq:smith-minors} becomes
\[
\mu_n \le \sum_{\abs{m} = n,\ i} \vp\bigl(\ip{m}{\lambda} - \lambda_i\bigr) - (M_n - 1)a.
\]
For rational $A$, Hadamard's inequality applied to a non-zero integer minor of order $r_n$ of $qL_n$, with $q$ prime to $p$, $qA$ integral and $H\ge1$ a bound for the absolute values of the entries of $qA$, gives the coarse bound $\mu_n\le r_n\log_p\bigl(\sqrt{r_n}(n+1)H\bigr)-(r_n-1)a$, which needs only the rank and also covers resonant and non-diagonalizable $A$.

\begin{proposition}[Semisimple linear part]\label{prop:smith-semisimple}
Let $L/\Qp$ be a finite extension with ring of integers $\mathcal O_L$, with $\vp$ normalized by $\vp(p)=1$, and suppose
$A=S\Lambda S^{-1}$, where $S\in\mathrm{GL}_d(L)$ and $\Lambda=\operatorname{diag}(\lambda_1,\dots,\lambda_d)$. Use the maximum-entry norm for matrices and put
\[
\kappa(S)=\log_p(\norm{S}\,\norm{S^{-1}}),\qquad
b_n=\max_{\substack{|m|=n,\ i\\ \langle m,\lambda\rangle\ne\lambda_i}}
\vp(\langle m,\lambda\rangle-\lambda_i).
\]
Then for every $j$, with $v_j(n)$ as in \eqref{eq:newton-hodge},
\begin{equation}\label{eq:smith-index}
\abs{s_j(n)-v_j(n)}\le(n+1)\kappa(S),
\end{equation}
where both sides are infinite for the same $j$. In particular, if $r_n>0$,
\begin{equation}\label{eq:smith-semisimple}
\mu_n\le b_n+(n+1)\kappa(S).
\end{equation}
If $S\in\mathrm{GL}_d(\mathcal O_L)$, the finite $s_j(n)$ are exactly the valuations of the non-zero $\langle m,\lambda\rangle-\lambda_i$, with multiplicity. In particular, in dimension one with $\lambda\ne0$,
\begin{equation}\label{eq:smith-one}
s_1(n)=\vp(\lambda)+\vp(n-1).
\end{equation}
\end{proposition}

\begin{proof}
On homogeneous vector polynomials use the coefficient Gauss norm and the change of coordinates $T_n h(z)=S^{-1}h(Sz)$. Then $T_nL_n=L_{\Lambda}|_{V_n}T_n$, and
\[
\norm{T_n}\le\norm{S^{-1}}\norm{S}^{n},\qquad
\norm{T_n^{-1}}\le\norm{S}\norm{S^{-1}}^{n}.
\]
The Smith exponents over $\mathcal O_L$, normalized by $\vp(p)=1$, are those over $\Zp$, since minors do not change. In the monomial eigenbasis $L_\Lambda|_{V_n}$ is diagonal with entries $\ip{m}{\lambda}-\lambda_i$, so its exponents are the $v_j(n)$. Hodge slopes are characterized by a min-max principle: $p^{-s_j(C)}$ is the smallest $\sigma$ such that every $j$-dimensional subspace $V$ contains $v\ne0$ with $\norm{Cv}\le\sigma\norm{v}$ \cite{kedlaya}. Hence, for invertible $B$ and $B'$, $s_j(BCB')\ge s_j(C)-\log_p\norm{B}-\log_p\norm{B'}$. Applied to $L_n=T_n^{-1}L_\Lambda T_n$ and to $L_\Lambda=T_nL_nT_n^{-1}$, this gives \eqref{eq:smith-index}, since $\log_p(\norm{T_n}\norm{T_n^{-1}})\le(n+1)\kappa(S)$. The two matrices have the same rank, so the largest finite $v_j(n)$ is $b_n$, and \eqref{eq:smith-semisimple} follows even in the presence of resonances. If $S$ and $S^{-1}$ are integral, then $\kappa(S)=0$, proving the equality assertion.
\end{proof}

\begin{corollary}[Algebraic eigenvalues]\label{cor:smith-algebraic}
In Proposition~\ref{prop:smith-semisimple}, suppose in addition that the $\lambda_i$ are algebraic over $\mathbb Q$. Let $k$ be a number field containing them, $d_k=[k:\mathbb Q]$, and choose a positive integer $q$ such that all $q\lambda_i$ are algebraic integers. Put
\[
H=\max\bigl(1,\max_{\tau:k\hookrightarrow\mathbb C,\,i}|\tau(q\lambda_i)|\bigr).
\]
Then, whenever $r_n>0$,
\begin{equation}\label{eq:smith-algebraic}
\mu_n\le(n+1)\kappa(S)+d_k\log_p((n+1)H)-\vp(q).
\end{equation}
In particular, an integral eigenbasis gives $\mu_n=O(\log(n+1))$. The algebraicity assumption here is over $\mathbb Q$, not merely over $\Qp$.
\end{corollary}

\begin{proof}
For a non-zero $\beta=\langle m,\lambda\rangle-\lambda_i$, the number $q\beta$ is an algebraic integer. Its non-zero integer norm has absolute value at most $((n+1)H)^{d_k}$. All its $p$-adic conjugates have non-negative valuation, so the valuation at the chosen embedding is at most the valuation of the norm. Hence
$\vp(\beta)\le d_k\log_p((n+1)H)-\vp(q)$. Apply \eqref{eq:smith-semisimple}.
\end{proof}

\begin{remark}[Why hypotheses on the coefficients matter]\label{rem:smith-liouville}
No logarithmic bound in $n$ holds for arbitrary $p$-adic coefficients, even for a diagonal $2\times2$ matrix. Choose integers $1\le b_1<b_2<\cdots$ with $b_{j+1}/p^{b_j}\to\infty$, put $\alpha=\sum_j p^{b_j}$ and $n_j=\sum_{k\le j}p^{b_k}$, and take $A=\operatorname{diag}(p^c,p^c\alpha)$, with $c$ satisfying \eqref{eq:hyp}. On $z_1^{n_j}\varepsilon_2$ the homological operator has eigenvalue $p^c(n_j-\alpha)$, of valuation $c+b_{j+1}$. Since the operator is diagonal in the integral monomial basis, this is a finite Smith exponent. Here $n_j\le p^{b_j+1}$, so the exponents grow faster than linearly along this subsequence, and choosing $b_{j+1}\ge f(p^{b_j+1})$ for a given increasing function $f$ makes them exceed $f(n_j)$. The number $\alpha$ is a $p$-adic Liouville number: $\vp(\alpha-n_j)=b_{j+1}$ grows faster than any multiple of $n_j$ (see \cite{kedlaya} for this notion). Such behavior explains the role of small-divisor hypotheses in full analytic normalization \cite{hermanyoccoz}, and it is the same obstruction as in the $p$-adic Fuchs theorem, where a regular singular linear system is normalized analytically when the differences of its exponents are non-Liouville, and Monsky's example shows that this hypothesis cannot be dropped \cite{kedlaya}. The threshold of Theorem~\ref{thm:threshold-convergence} avoids it by keeping such divisors in $R$, as Liouville partitions keep exponents with Liouville differences in one block for linear connections \cite[Definition~3.4.4]{kedlaya15}. The integer and number-field bounds above do not apply to this choice of coefficients.
\end{remark}

For a linear part $\lambda I + J$ with $\lambda\ne0$ and $J$ nilpotent, $L_n$ is $L_{\lambda I}$ plus a nilpotent operator, so $L_n^{-1}$ is a finite Neumann series and its exponents can be bounded explicitly. In the plane they can even be read off from weighted paths. I do not need such bounds below, since the threshold keeps whatever is too small.

\subsection{Convergence and tail bounds}\label{sec:normal-tail}

If only directions with $s_i(n)\le\sigma$ are removed, each degree loses at most $\sigma$ digits, whatever the bounds above give, and the normalization converges with no Diophantine assumption on $A$.

\begin{theorem}[Normalization with a fixed threshold]\label{thm:threshold-convergence}
Let $\tilde f(z)=Az+\sum_{n\ge2}\tilde f_n(z)$ satisfy
$\vp(A)\ge c$ and $\vp(\tilde f_n)\ge c+r(n-1)$, where $c>1/(p-1)$ and $r\ge1$ are integers. Fix an integer $\sigma\ge0$ and use the Smith construction above at every degree, removing only directions with $s_i(n)\le\sigma$. Set
\[
b=\max(c,\sigma),\qquad \eta=r+c-b=r-\max(\sigma-c,0).
\]
If $\eta\ge1$, all divisions are integral, and for every $n\ge2$,
\begin{equation}\label{eq:normal-majorant}
\vp(G_n)\ge b+\eta(n-1),\qquad
\vp(\Phi_n)\ge\eta(n-1),\qquad
\vp(R_n)\ge b+\eta(n-1).
\end{equation}
The series $\Phi$ and $g_{\mathrm{nf}}$ converge in the Tate algebra on $\Zp^d$, satisfy \eqref{eq:conj}, and $\Phi$ is an analytic bijection of $\Zp^d$ with analytic inverse. For $M\ge1$ their tails satisfy
\begin{equation}\label{eq:normal-tail}
\vp(\Phi-\Phi_{\le M})\ge\eta M,\qquad
\vp(g_{\mathrm{nf}}-g_{\mathrm{nf},\le M})\ge b+\eta M.
\end{equation}
Thus $r\ge1+\max(\sigma-c,0)$ is a sufficient radius exponent for this partial normalization, including non-semisimple and resonant linear parts.
\end{theorem}

\begin{proof}
Assume the bounds in \eqref{eq:normal-majorant} through degree $n-1$. In a degree-$n$ contribution from $\tilde f_l(\Phi)$, the $l$ factors from $\Phi$ have degrees summing to $n$. The linear factor has valuation zero, so the coefficient valuation of this contribution is at least
\[
c+r(l-1)+\eta(n-l)
=\eta(n-1)+c+(b-c)(l-1)
\ge\eta(n-1)+b.
\]
In \eqref{eq:homological}, the other contributions to $G_n$ are $D\Phi_j\cdot R_k$ with $j+k=n+1$ and $2\le j,k<n$. Differentiation preserves a lower bound on coefficient valuations, and the induction gives $\eta(j-1)+b+\eta(k-1)=b+\eta(n-1)$ for these terms. This also starts the induction at $n=2$, when there are no such products.

The matrices $P,Q,P^{-1}$ are integral, so the only loss in forming $\Phi_n$ or $R_n$ comes from the divisions by $p^{s_i}$ in the removed directions. There the loss is at most $\sigma\le b$, giving $\vp(\Phi_n)\ge\eta(n-1)$, while the retained directions give $\vp(R_n)\ge b+\eta(n-1)$. Since $\eta\ge1$, these bounds give convergence and \eqref{eq:normal-tail}. The conjugacy identity then holds in the Tate algebra, and $b\ge c$ puts $g_{\mathrm{nf}}$ in $p^c\Zp\langle z\rangle^d$.

To invert $\Phi=\mathrm{id}+H$, solve $z=w-H(z)$. The coefficients of $H$ are divisible by $p$, which makes this a contraction both pointwise on $\Zp^d$ and on the unit ball of the Tate algebra in $w$. Starting with $z=w$ therefore gives an inverse in the Tate algebra.
\end{proof}

The radius is conservative, since it charges every removed direction with the full $\sigma$ digits. Letting the loss depend on the degree refines it (in dimension one it recovers the bound for $h$ at the start of this section), but the certificate of Section~\ref{sec:normal-certificate} makes such refinements unnecessary in practice: it reads the actual valuations of the computed coefficients.

\begin{lemma}[Stability of the flow under coefficient truncation]\label{lem:flow-stability}
Let $g,h\in p^c\Zp\langle z\rangle^d$, with $c>1/(p-1)$, and suppose $\vp(g-h)\ge L$. Then their flows agree modulo $p^L$ for all $t\in\Zp$ and $z\in\Zp^d$.
\end{lemma}

\begin{proof}
Write $D_g=g\cdot\nabla$ and $D_h=h\cdot\nabla$. The difference of the two Lie series is controlled term by term by the identity
\[
D_g^j-D_h^j=\sum_{a=0}^{j-1}D_g^{j-1-a}(D_g-D_h)D_h^a
\]
because each summand contains one factor that raises coefficient valuations by at least $L$ and $j-1$ factors that raise them by at least $c$. Applied to $z$ and divided by $j!$, it has valuation at least $L+(j-1)c-\vp(j!)\ge L$. Both series converge, so the same bound holds for their difference.
\end{proof}

\begin{corollary}[Certified truncation]\label{cor:normal-cut}
Under Theorem~\ref{thm:threshold-convergence}, let $K\ge1$ and $M\ge1$ satisfy $\eta M\ge K$. Let $\psi_t^{(M)}$ be the flow of the polynomial field $g_{\mathrm{nf},\le M}$. Then the flow $\tilde\varphi_t$ of $\tilde f$ satisfies
\begin{equation}\label{eq:normal-certified-flow}
\tilde\varphi_t(z)\equiv
\Phi_{\le M}\!\left(\psi_t^{(M)}(\Phi_{\le M}^{-1}(z))\right)
\pmod{p^K}
\end{equation}
for all $t\in\Zp$ and $z\in\Zp^d$. The inverse on the right is computed by contraction, and all polynomials may be reduced coefficientwise modulo $p^K$ after their coefficients have been determined to that precision.
\end{corollary}

\begin{proof}
The exact analytic conjugacy gives $\tilde\varphi_t=\Phi\circ\psi_t\circ\Phi^{-1}$. The bound \eqref{eq:normal-tail} gives $\Phi\equiv\Phi_{\le M}\pmod{p^K}$ and $g_{\mathrm{nf}}\equiv g_{\mathrm{nf},\le M}\pmod{p^K}$. Both coordinate maps are isometries, being the identity plus a strict contraction. Their inverses therefore agree modulo $p^K$ as well. By Lemma~\ref{lem:flow-stability}, $\psi_t$ and $\psi_t^{(M)}$ agree modulo $p^K$. All maps in the composition preserve congruences, proving \eqref{eq:normal-certified-flow}. Returning to the original variables uses $x=x^*+p^r z$.
\end{proof}

\subsection{A certificate for the computed normal form}\label{sec:normal-certificate}

The exact construction gives a tail bound for its own coefficients. In a calculation modulo $p^W$, a quotient $\gamma_i/p^{s_i}$ is determined modulo $p^{W-s_i}$. The lift chosen in one degree enters the equations in later degrees, and the Smith matrices are themselves computed modulo $p^W$. I check the computed polynomials by evaluating the defect in \eqref{eq:conj}.

\begin{lemma}[Conjugacy defect]\label{lem:defect}
Let $c>1/(p-1)$ and $K\ge c$ be integers, $\tilde f\in p^c\Zp\langle z\rangle^d$, and let $\Phi=\mathrm{id}+H$ and $g_{\mathrm{nf}}$ be polynomial maps with $H\in p\Zp[z]^d$ and $g_{\mathrm{nf}}\in p^c\Zp[z]^d$. If $\mathcal D=D\Phi\cdot g_{\mathrm{nf}}-\tilde f\circ\Phi$ satisfies $\vp(\mathcal D)\ge K$, then $\tilde\varphi_t(\Phi(w))\equiv\Phi(\psi_t(w))\pmod{p^K}$ for all $t\in\Zp$ and $w\in\Zp^d$, where $\psi_t$ is the flow of $g_{\mathrm{nf}}$.
\end{lemma}

\begin{proof}
Transporting $g_{\mathrm{nf}}$ by $\Phi$ gives the field $g=\tilde f+\mathcal D\circ\Phi^{-1}$, since $g\circ\Phi=D\Phi\cdot g_{\mathrm{nf}}$. The contraction argument of Theorem~\ref{thm:threshold-convergence} gives $\Phi^{-1}\in\mathrm{id}+p\Zp\langle z\rangle^d$ and shows that $\Phi$ is an isometry. Substitution by $\Phi^{-1}$ therefore preserves the coefficient valuation bound on the defect, giving $g\in p^c\Zp\langle z\rangle^d$ and $\vp(g-\tilde f)\ge K$.

By the chain rule, $t\mapsto\Phi(\psi_t(w))$ solves $x'=g(x)$ with $x(0)=\Phi(w)$. Uniqueness among power series solutions identifies it with the flow of $g$ (Section~\ref{sec:flow}), and Lemma~\ref{lem:flow-stability} gives the comparison with the flow of $\tilde f$.
\end{proof}

\begin{lemma}[The defect above the truncation degree]\label{lem:defect-tail}
Let $\Phi$ and $g_{\mathrm{nf}}$ be polynomial maps of degree at most $M$, with $\Phi-\mathrm{id}\in p\Zp[z]^d$. Suppose there are $b\ge0$ and $\eta>0$ with
\[
\vp(\Phi_n)\ge\eta(n-1),\qquad \vp(g_{\mathrm{nf},n})\ge b+\eta(n-1)\quad(2\le n\le M),\qquad \vp(\tilde f_n)\ge b+\eta(n-1)\quad(n\ge2).
\]
Then the homogeneous parts of $\mathcal D$ satisfy $\vp(\mathcal D_n)\ge b+\eta(n-1)$ for every $n>M$.
\end{lemma}

\begin{proof}
In degree $n>M$, the only possible terms in $D\Phi\cdot g_{\mathrm{nf}}$ are $D\Phi_j\cdot g_{\mathrm{nf},k}$ with $j+k=n+1$ and $2\le j,k\le M$. They have valuation at least $\eta(j-1)+b+\eta(k-1)=b+\eta(n-1)$. The terms of $\tilde f$ of degree at most one contribute nothing in these degrees.

For the contribution of $\tilde f_l\circ\Phi$ with $l\ge2$, a product in degree $n$ contains a coefficient of $\tilde f_l$ and $l$ coefficients of $\Phi$ whose degrees $k_1,\dots,k_l$ sum to $n$. Each factor from $\Phi$ has valuation at least $\eta(k_i-1)$, including $k_i=0,1$ by the condition on $\Phi-\mathrm{id}$. The resulting valuation is at least $b+\eta(l-1)+\eta(n-l)=b+\eta(n-1)$.
\end{proof}

\begin{theorem}[Certified normal form]\label{thm:certificate}
Fix an integer $K\ge1$, with $K\ge2$ if $p=2$. For $\tilde x\in\Zp^d$ and $r\ge1$, put $\tilde f(z)=p^{-r}f(\tilde x+p^rz)$. Suppose $\Phi=\mathrm{id}+H$ and $g_{\mathrm{nf}}$ are polynomial maps of degree at most $M$ such that
\begin{enumerate}[leftmargin=*,label=(\roman*)]
\item $H\in p\Zp[z]^d$, and $g_{\mathrm{nf}}$ and $\tilde f$ have coefficients divisible by $p$ (by $4$ if $p=2$);
\item $\mathcal D_n\equiv0\pmod{p^K}$ for $0\le n\le M$;
\item the bounds of Lemma~\ref{lem:defect-tail} hold with $b+\eta M\ge K$, where coefficients of degree $\ge2$ that vanish modulo $p^K$ are ignored.
\end{enumerate}
Then for all $t\in\Zp$ and all $x$ with $\abs{x-\tilde x}\le p^{-r}$,
\begin{equation}\label{eq:certified}
\varphi_t(x)\equiv\tilde x+p^r\,\Phi\bigl(\psi_t\bigl(\Phi^{-1}\bigl((x-\tilde x)/p^r\bigr)\bigr)\bigr)\pmod{p^{K+r}}.
\end{equation}
\end{theorem}

\begin{proof}
Delete the coefficients of degree $\ge2$ that vanish modulo $p^K$, and write $\hat\Phi,\widehat g_{\mathrm{nf}},\hat f$ for the resulting polynomials or series. Their defect $\widehat{\mathcal D}$ differs from $\mathcal D$ by terms containing coefficients of valuation at least $K$. By (ii), $\vp(\widehat{\mathcal D}_n)\ge K$ through degree $M$, and the literal bounds in (iii) give the same inequality in higher degrees by Lemma~\ref{lem:defect-tail}.

This gives $\vp(\mathcal D)\ge K$ for the original polynomials. Use $c=1$ at odd primes and $c=2$ at $p=2$ in Lemma~\ref{lem:defect}. Applied at $w=\Phi^{-1}(z)$, it gives the congruence in rescaled coordinates. Multiplying by $p^r$ gives \eqref{eq:certified}.
\end{proof}

Conditions (i)--(iii) involve only finitely many coefficients and are checked in exact integer arithmetic. The center $\tilde x$ need not be an exact equilibrium: an error in it appears as the constant term of $\tilde f$, which (ii) requires to vanish modulo $p^K$. For a power series field the coefficients beyond the cut enter only through a declared lower bound for their valuations.

For the computation take $W=K+\sigma+2$. In each degree the computed matrices satisfy $PL_nQ\equiv\operatorname{diag}(p^{s_i})\pmod{p^W}$. I form $G_n$ from the lower-degree polynomials already computed, divide the residues in directions with $s_i\le\sigma$ by $p^{s_i}$, and retain the other directions in $R_n$. The resulting equation $L_A\Phi_n+R_n\equiv G_n\pmod{p^W}$ gives $\mathcal D_n\equiv0\pmod{p^W}$ for $1\le n\le M$.

After each degree I use the computed valuations to look for $b\in[0,K]$ and $\eta>0$ satisfying Lemma~\ref{lem:defect-tail}. The construction stops at the first degree $M$ for which $b+\eta M\ge K$. To see when this is guaranteed, use the induction of Theorem~\ref{thm:threshold-convergence}. Taking canonical residues modulo $p^W$ preserves lower bounds on valuations. Dividing a nonzero residue by $p^{s_i}$ lowers its valuation by exactly $s_i$. The proof therefore gives the same bounds on the computed polynomials. If $f(\tilde x)\equiv0\pmod{p^{W+r}}$ and $r=1+\max(\sigma-c,0)$, the bounds hold with $b=\max(c,\sigma)$ and $\eta=1$. If this value of $b$ exceeds $K$, it can be lowered to $K$. Condition (i) of Theorem~\ref{thm:certificate} holds as well, and the construction succeeds at degree $M\le\max(K,2)$, since it starts at degree $2$. The actual coefficient valuations often allow a smaller $r$.

For example, take $p=3$ and $f(x)=3x+3x^2$ at $0$. With $r=1$ the rescaled field is $\tilde f(z)=3z+9z^2$. In degree two the homological operator is multiplication by $3$, so a threshold $\sigma=1$ gives $\Phi(z)=z+3z^2$ and $g_{\mathrm{nf}}(z)=3z$. For $K=3$ and $M=2$ the defect is $-54z^3-81z^4$, with no terms through degree two. The quadratic coefficients of $\Phi$ and $\tilde f$ have valuations $1$ and $2$, respectively, so $b=\eta=1$ gives $b+\eta M=3$. The certificate compares the flows modulo $3^4$ on $3\mathbb Z_3$.

\subsection{Block structure of the homological operator}\label{sec:normal-blocks}

A dense Smith form of $L_n$ costs $O(M_n^3)$ operations on $W$-digit integers, with $M_n=d\binom{n+d-1}{d-1}$. This is what limits the degree and the dimension. The operator is only as large as the coupling of the coordinates by $A$ requires.

\begin{proposition}[Residual splitting]\label{prop:residual-splitting}
Let $A=p^aA'$ with $A'$ integral and $A'\not\equiv0\pmod p$, and let $\bar\chi=G_1\cdots G_B$ be a factorization of the characteristic polynomial of $A'$ modulo $p$ into pairwise coprime factors. There is $T\in\mathrm{GL}_d(\Zp)$ such that $T^{-1}AT=\operatorname{diag}(A_1,\dots,A_B)$, where $A_b$ has size $d_b=\deg G_b$ and $p^{-a}A_b$ has residual characteristic polynomial $G_b$. In the coordinates $w=T^{-1}z$, for $\kappa=(n_1,\dots,n_B)$ with $\abs{\kappa}=n$ and $1\le b\le B$, let $V_{n,\kappa,b}$ be spanned by the $w^e\varepsilon_i$ whose exponents have total degree $n_\beta$ in block $\beta$ and with $i$ in block $b$. These spaces are invariant under $L_n$, of dimension $d_b\prod_\beta\binom{n_\beta+d_\beta-1}{d_\beta-1}$, and the Smith exponents of $L_n$ are the union of those of its restrictions. If $G_\beta=(x-\alpha_\beta)^{d_\beta}$ for all $\beta$ and $\sum_\beta n_\beta\alpha_\beta\not\equiv\alpha_b\pmod p$, every exponent on $V_{n,\kappa,b}$ equals $a$.
\end{proposition}

\begin{proof}
By Hensel's lemma the factorization lifts to $\Zp[x]$, and the corresponding idempotents of $\Zp[A']$ split $\Zp^d$ into $A$-invariant summands of ranks $d_b$. The matrix $T$ maps the coordinate blocks onto them. In the new coordinates $A$ is block diagonal, so the derivation part of $L_A$ keeps the degree of a monomial in each block and the term $A\Phi$ keeps the block of $i$. Changing coordinates by $T$ is unimodular on $V_n$, so the exponents are those of the block diagonal operator. In the last case $L_{A'}$ has the single eigenvalue $\sum_\beta n_\beta\alpha_\beta-\alpha_b\ne0$ modulo $p$ on $V_{n,\kappa,b}$, hence is invertible there over $\Zp$, and $L_A=p^aL_{A'}$.
\end{proof}

So for a residually split $A$ the directions with exponent above $a$ are confined to the pieces with $\sum_\beta n_\beta\alpha_\beta\equiv\alpha_b$, a condition on the residues of the block degrees modulo $p$. For $d=1$ it is $p\mid n-1$. One need not factor $\bar\chi$: for $B=A'-\alpha$ with $\alpha\in\mathbb F_p$, or $B=A'^{p^k}-A'$ for the eigenvalues of degree dividing $k$ over $\mathbb F_p$, the limit $\lim_j(B^{e_0})^{p^j}$ with $e_0=p^{\lceil\log_pd\rceil}\operatorname{lcm}_{i\le d}(p^i-1)$ is the idempotent onto the part on which $B$ is invertible modulo $p$, and $T$ is assembled from columns of these projectors. The recursion runs piece by piece in the coordinates $w$, and the certificate of Section~\ref{sec:normal-certificate}, which refers to the original coordinates, does not depend on $T$.

On the pieces the Smith exponents agree with those of the whole $V_n$, and the cost of the decompositions grows with the largest piece instead of $M_n$.

The certificate concerns the flow modulo $p^K$ on the chosen ball. For example, in $(px,\,2py+p^Kx^2)$ the resonant term is invisible at this precision and survives in the exact field. Thus a computed linear normal form can pass the check while the exact field has a resonant term. For an analytic input, the certificate also depends on the supplied valuation bound for its tail.

\section{A \texorpdfstring{$p$}{p}-adic Eliasson--Vey theorem}\label{sec:ev}

Putting the resonant terms into a function of the quadratic parts makes a Newton iteration possible for an integrable system. Its divisors in symplectic eigencoordinates are integers, so their valuation grows at most logarithmically with degree and the losses are summable as the degrees double. This yields an explicit radius and answers the rank-$0$ part of Crespo and Pelayo's question about Williamson types \cite[Question~5.14]{crespopelayoww}.

On $\Qp^{2n}$ with its standard symplectic form, consider an analytic system $\mathcal H=(H_1,\dots,H_n)$ with $\{H_i,H_j\}=0$ near $0$. Assume $\mathcal H(0)=0$ and $d\mathcal H(0)=0$, so the origin is a critical point of rank $0$. The point is \emph{non-degenerate} if the quadratic parts of the $H_i$, viewed in $\mathfrak{sp}(2n,\Qp)$, span an $n$-dimensional Cartan subalgebra $\mathfrak t$ of commuting semisimple elements. Choose a basis $q=(q_1,\dots,q_n)$ of $\mathfrak t$, for instance a Williamson basis in the sense of \cite{crespopelayoww}. After a linear change on the target, I take $\mathcal H$ to have quadratic part $q$.

\begin{theorem}\label{thm:ev}
There are an analytic symplectic map $\Phi$ from a ball around $0$ onto a neighborhood of $0$, with $\Phi(0)=0$, and an analytic map $\mathcal G\colon(\Qp^n,0)\to(\Qp^n,0)$ with $d\mathcal G(0)=\mathrm{id}$, both with coefficients in $\Qp$, such that
\[
\mathcal H\circ\Phi=\mathcal G(q_1,\dots,q_n).
\]
In the eigencoordinates of Lemma~\ref{lem:ev-coords} the ball has radius $p^{-\rho}$, so in the original coordinates it contains the ball of radius $p^{-\rho}\norm{S^{-1}}^{-1}$, where $S^{-1}$ is the linear map to the eigencoordinates. Here $\rho$ is any integer with $\rho>3r+\mu+L_p+\frac1{p-1}$, $r$ and $\mu$ measure $\mathcal H$ as below, and
\[
L_p=\sum_{k\ge0}\frac{\lfloor\log_p(2^{k+1}+1)\rfloor+\frac1{p-1}}{2^k},\qquad L_2=6,\quad L_3\approx3.68,\quad L_5\approx1.66,\quad L_7\approx0.905.
\]
\end{theorem}

\begin{lemma}[Eigencoordinates]\label{lem:ev-coords}
There are a finite Galois extension $L/\Qp$, with valuation normalized by $\vp(p)=1$, and $L$-linear symplectic coordinates $(z_1,\dots,z_n,w_1,\dots,w_n)$, $\{z_i,w_j\}=\delta_{ij}$, such that with $\zeta_i=z_iw_i$ one has $q=\Lambda\zeta$ for some $\Lambda\in\mathrm{GL}_n(L)$ and
\[
\{\zeta_i,z^\alpha w^\beta\}=(\beta_i-\alpha_i)\,z^\alpha w^\beta .
\]
\end{lemma}

\begin{proof}
The Cartan subalgebra splits over a finite extension, and I take $L$ to be its Galois closure. Then $L^{2n}$ is a sum of weight lines with pairwise distinct non-zero weights $\pm\varepsilon_1,\dots,\pm\varepsilon_n$. Invariance of the symplectic form $\omega$ gives $(\varepsilon+\varepsilon')(T)\,\omega(x,y)=0$ for $x,y$ of weights $\varepsilon,\varepsilon'$ and $T\in\mathfrak t$, so $\omega$ pairs the line of $\varepsilon_i$ only with that of $-\varepsilon_i$, non-degenerately. Choosing $e_i,e_i^\prime$ with $\omega(e_i,e_i^\prime)=1$ gives a symplectic basis, and the dual coordinates satisfy the brackets above. The quadratic Hamiltonian of $T\in\mathfrak t$ is $\sum_i\varepsilon_i(T)z_iw_i$, and $T\mapsto(\varepsilon_i(T))_i$ is a bijection onto $L^n$, which gives $\Lambda$.
\end{proof}

So the joint homological operator $W\mapsto(\{\zeta_i,W\})_i$ is diagonal on monomials, with the integer vector $\delta=\beta-\alpha$ as eigenvalue, and its kernel is $L[\zeta]=L[q]$. A non-zero integer of absolute value at most $D$ is divisible by at most $p^{\lfloor\log_pD\rfloor}$, so for a monomial of degree $D$ with $\delta\ne0$,
\begin{equation}\label{eq:ev-divisor}
\min\{\vp(\delta_i):\delta_i\ne0\}\le\lfloor\log_pD\rfloor .
\end{equation}
These integers have complex absolute value at least $1$, which accounts for the absence of small divisors in Vey's theorem. The $p$-adic cost grows at most logarithmically with the degree.

To keep track of coefficient losses, give a homogeneous term $f_D$ of degree $D$ in $(z,w)$ weight $D-2$. For $a>0$, let $\mathcal A_a$ consist of series $f=\sum_{D\ge2}f_D$ over $L$ with $\vp(f_D)\ge a(D-2)$. They converge on the closed unit polydisk. In the variables $u\in L^n$, write $h=\sum_\gamma h_\gamma u^\gamma\in\mathcal U_a$ when $\vp(h_\gamma)\ge a(2\abs\gamma-2)$. Each $\zeta_i$ has degree two, so $h(\zeta)\in\mathcal A_a$ if $h(0)=0$. Write $\operatorname{ad}_W f=\{f,W\}$.

\begin{lemma}\label{lem:ev-classes}
Let $a>0$.
\begin{enumerate}[label=(\alph*),leftmargin=*]
\item $\{\mathcal A_a,\mathcal A_a\}\subset\mathcal A_a$.
\item If $h\in\mathcal U_a$, $h(0)=0$, and $y_1,\dots,y_n\in\mathcal A_a$ have no terms of degree below $2$, then $h(y)\in\mathcal A_a$.
\item If $h\in\mathcal U_a$ has no terms of degree below $2$, then $u\mapsto u+h(u)$ has an inverse $u\mapsto u-\tilde h(u)$ with $\tilde h\in\mathcal U_a$.
\item If $W\in\mathcal A_b$ has no terms of degree below $N\ge3$, $f\in\mathcal A_b$ and $k\ge1$, then $\operatorname{ad}_W^kf/k!$ has no terms of weight below $k(N-2)$ and lies in $\mathcal A_{b-1/((p-1)(N-2))}$.
\end{enumerate}
\end{lemma}

\begin{proof}
The degree of a bracket $\{f_D,g_{D^{\prime}}\}$ is $D+D^\prime-2$, so its weight is the sum of the two weights. Differentiation only multiplies coefficients by integers, and the bracket has valuation at least $a(D-2)+a(D^\prime-2)$. This gives (a).

In (b), the product of $k$ homogeneous factors of total degree $D=\sum D_i$ has valuation at least $a(D-2k)$. The coefficient $h_\gamma$ with $\abs\gamma=k$ supplies the extra $a(2k-2)$. The same count, with each $u_i$ assigned degree two and weight zero, applies to substitutions in the equation $U=y-h(U)$. Iterating from $U=y$ stays in $\mathcal U_a$ and determines the inverse coefficient by coefficient, giving (c).

For (d), apply (a) at every bracket. The weight is at least $k(N-2)$, and dividing by $k!$ costs at most $(k-1)/(p-1)\le\text{weight}/((p-1)(N-2))$, giving the stated bound.
\end{proof}

If moreover $b(N-2)>1/(p-1)$, then $\exp(\operatorname{ad}_W)f=\sum_k\operatorname{ad}_W^kf/k!$ converges for $f\in\mathcal A_b$ and for the coordinate functions: for a coordinate, $\operatorname{ad}_W^k$ gives parts of degree $D\ge1+k(N-2)$ and valuation at least $b(D-1)$, and dividing by $k!$ leaves at least $\bigl(b-\frac1{(p-1)(N-2)}\bigr)(D-1)$. So $\exp(\operatorname{ad}_W)$ is an automorphism of the Poisson algebra, with inverse $\exp(-\operatorname{ad}_W)$, and $\phi_W=(\exp(\operatorname{ad}_W)z,\exp(\operatorname{ad}_W)w)$ is an analytic symplectic automorphism of the unit polydisk with $\exp(\operatorname{ad}_W)f=f\circ\phi_W$. The Lie series constructs the time-one map directly.

\begin{proposition}[The Newton step]\label{prop:ev-step}
Let $\mathcal H=\zeta+R$ with $\{H_i,H_j\}=0$ and $R\in\mathcal A_a$ without terms of degree below $N\ge3$. Put $s_N=\lfloor\log_p(2N-3)\rfloor$ and $b=a-s_N/(N-2)$, and assume $b(N-2)>1/(p-1)$. There are $W\in\mathcal A_b$, homogeneous of degrees in $[N,2N-3]$, and $h\in\mathcal U_a$ such that, with $\mathcal G(u)=u+h(u)$,
\[
\mathcal H'=\mathcal G^{-1}\circ \mathcal H\circ\phi_W=\zeta+R',\qquad R'\in\mathcal A_{a'},\quad a'=a-\frac{s_N+\frac1{p-1}}{N-2},
\]
$R'$ has no terms of degree below $2N-2$, and $\{H'_i,H'_j\}=0$.
\end{proposition}

\begin{proof}
The commutation of the $H_i$ gives $\{\zeta_i,R_j\}-\{\zeta_j,R_i\}=-\{R_i,R_j\}$. Since the right side starts in degree $2N-2$, a monomial $m=z^\alpha w^\beta$ below that degree, with $\delta=\beta-\alpha$, satisfies
\begin{equation}\label{eq:ev-compat}
\delta_iR_{j,m}=\delta_jR_{i,m}\qquad\text{for all }i,j.
\end{equation}

Only degrees $N$ through $2N-3$ enter the correction. If $\delta\ne0$, take $W_m=-R_{i,m}/\delta_i$ with $\delta_i\ne0$. The relation \eqref{eq:ev-compat} makes this independent of $i$ and forces $R_{j,m}=0$ whenever $\delta_j=0$. Thus $\{\zeta_i,W\}=-R_i^{\mathrm{nr}}$, where $R^{\mathrm{nr}}$ consists of those low-degree terms with $\delta\ne0$. Choosing the least valuation among the nonzero $\delta_i$ gives $\vp(W_m)\ge a(D-2)-s_N\ge b(D-2)$ for $D=\deg m$, by \eqref{eq:ev-divisor}. In particular $W\in\mathcal A_b$.

The Lie series removes $R_i^{\mathrm{nr}}$ from the low degrees of $H_i\circ\phi_W$, leaving $\zeta_i+R_i^{\mathrm{res}}+\mathcal R_i$. Here $R^{\mathrm{res}}$ consists of the terms in degrees $N$ through $2N-3$ with $\delta=0$, and
\[
\mathcal R_i=R_i^{\ge2N-2}+\sum_{k\ge1}\frac{\operatorname{ad}_W^kR_i}{k!}+\sum_{k\ge2}\frac{\operatorname{ad}_W^{k-1}(-R^{\mathrm{nr}}_i)}{k!}.
\]
Bracketing with $W$ adds at least $N-2$ to the degree, so $\mathcal R_i$ starts in degree $2N-2$. In weight $w$, the factors from $W$ cost at most $(a-b)w=s_Nw/(N-2)$ against the original slope $a$. Every term divided by $k!$ has weight at least $k(N-2)$, so the factorial costs at most $w/((p-1)(N-2))$. This gives the bound $a'w$ and hence $\mathcal R_i\in\mathcal A_{a'}$.

Write $R_i^{\mathrm{res}}=h_i(\zeta)$, taking $h_i(u)=\sum_\gamma R_{i,z^\gamma w^\gamma}u^\gamma$ over $2\abs\gamma\in[N,2N-3]$. These terms have $\abs\gamma\ge2$ and $h\in\mathcal U_a$, so $\mathcal H\circ\phi_W=\mathcal G(\zeta)+\mathcal R$ with $\mathcal G(u)=u+h(u)$. Lemma~\ref{lem:ev-classes}(c) gives $\mathcal G^{-1}=\mathrm{id}-\tilde h$, and then
\[
\mathcal H'-\zeta=\mathcal R-\bigl(\tilde h(\mathcal G(\zeta)+\mathcal R)-\tilde h(\mathcal G(\zeta))\bigr),
\]
where each term of the difference contains a factor of $\mathcal R$. This preserves the minimum degree $2N-2$, and Lemma~\ref{lem:ev-classes}(b) gives $\mathcal H'-\zeta\in\mathcal A_{a'}$. The change $\phi_W$ is symplectic and $\mathcal G^{-1}$ is a change on the target, so the components of $\mathcal H'$ still Poisson-commute.
\end{proof}

\begin{proof}[Proof of Theorem~\ref{thm:ev}]
In the coordinates of Lemma~\ref{lem:ev-coords}, $\Lambda^{-1}\mathcal H=\zeta+R$ with $R=\sum_{D\ge3}R_D$, and since $\mathcal H$ is analytic there are $r,\mu$ with $\vp(R_D)\ge-\mu-rD$ and $\mu+2r\ge0$. They include the non-integrality of the linear changes of coordinates. The rescaling $(z,w)\mapsto(p^\rho z,p^\rho w)$, $\mathcal H\mapsto p^{-2\rho}\mathcal H(p^\rho\,\cdot)$, preserves the bracket and $\zeta$ and multiplies $R_D$ by $p^{\rho(D-2)}$. Since $\rho(D-2)-\mu-rD-(\rho-3r-\mu)(D-2)=(D-3)(2r+\mu)\ge0$, the rescaled $R$ lies in $\mathcal A_{a_0}$ with $a_0=\rho-3r-\mu>L_p+\frac1{p-1}$.

Start with $\mathcal H_0=\Lambda^{-1}\mathcal H$ in the rescaled coordinates and put $N_k=2^k+2$, so that $N_{k+1}=2N_k-2$ and $N_k-2=2^k$, and $a_{k+1}=a_k-(s_{N_k}+\frac1{p-1})/2^k$ with $s_{N_k}=\lfloor\log_p(2^{k+1}+1)\rfloor$. Then $a_k\ge a_0-L_p>\frac1{p-1}$ for all $k$, and $b_k\ge a_{k+1}$, so the hypothesis of Proposition~\ref{prop:ev-step} holds at every step. It gives $W_k$, $\mathcal G_k$ and $\mathcal H_{k+1}=\mathcal G_k^{-1}\circ \mathcal H_k\circ\phi_{W_k}=\zeta+R_{k+1}$ with $R_{k+1}\in\mathcal A_{a_{k+1}}\subset\mathcal A_{a_0-L_p}$ and no terms of degree below $N_{k+1}$, so
\[
\mathcal H\circ\phi_{W_0}\circ\cdots\circ\phi_{W_{k}}=\Lambda\,\mathcal G_0\circ\cdots\circ\mathcal G_{k}(\zeta+R_{k+1}).
\]
In the spatial change, a homogeneous part of degree $D$ of $\phi_{W_k}-\mathrm{id}$ has $D\ge N_k-1$ and valuation at least $a_{k+1}(D-1)$ by the coordinate estimate for the Lie series. In the target change, a term $u^\gamma$ of $\mathcal G_k-\mathrm{id}$ has $2\abs\gamma\ge N_k$ and valuation at least $a_k(2\abs\gamma-2)$. Both estimates give the lower bound $(a_0-L_p)(N_k-2)$, positive from the first step and tending to infinity. The inverse changes satisfy the same estimates, by using $-W_k$ for the spatial inverse and Lemma~\ref{lem:ev-classes}(c) for the target inverse. The changes and their inverses are isometries, so their compositions converge to analytic automorphisms of the corresponding unit polydisks. The spatial limit $\Phi$ is symplectic.

Since $R_k\to0$, the equation above gives the desired identity in the rescaled coordinates with target map $\Lambda\mathcal G(\Lambda^{-1}\,\cdot)$, where $\mathcal G$ is the limit of the target compositions. Undoing the rescaling sends $\Phi$ to $p^\rho\Phi(p^{-\rho}\,\cdot)$, which is symplectic on the ball of radius $p^{-\rho}$. By $q(p^{-\rho}x)=p^{-2\rho}q(x)$, the corresponding target map is $u\mapsto p^{2\rho}\Lambda\mathcal G(p^{-2\rho}\Lambda^{-1}u)$ and still has derivative $\mathrm{id}$ at $0$.

To descend the construction, work in the original spatial coordinates and use the operators $\operatorname{ad}_{q_j}$ on each homogeneous degree. Their joint kernel and the sum of their images are complementary, as can be seen on the monomials in eigencoordinates. Both spaces are stable under $\operatorname{Gal}(L/\Qp)$ because the $q_j$ are defined over $\Qp$. Given an input over $\Qp$, the nonresonant equations have a unique solution $W_k$ in the sum of the images, and the resonant projection lies in the joint kernel. Both are therefore Galois invariant. Each component of the resonant projection lies in $L[q]$ and has a unique expression $P(q)$ with $P\in L[u]$, since the $q_j$ are algebraically independent. Galois invariance of $P(q)$ fixes every coefficient of $P$, and the fixed field is $\Qp$. Thus the spatial change and the correction on the original target, $\Lambda\mathcal G_k(\Lambda^{-1}\,\cdot)$, are defined over $\Qp$ at each step. The same holds for the limits.
\end{proof}

At a point of rank $0$ the Williamson normal form $(g_1,\dots,g_n)$ of \cite{crespopelayoww} is a basis of $\mathfrak t$, so Theorem~\ref{thm:ev} with $q=g$ gives exactly the local model asked for in \cite[Question~5.14]{crespopelayoww}, with no condition on the type. Non-degeneracy is used here in the Cartan sense above, which is the definition used in the real theory. Definition~5.2 of \cite{crespopelayoww} asks for linearly independent Hessians and a combination $\Omega^{-1}\sum_ia_i\,d^2H_i(0)$ with $n$ distinct eigenvalues, and Remark~5.3 there says that it is meant to agree with the Cartan definition. I read it with $2n$ distinct eigenvalues, as in \cite[Definition~3.2]{crespopelayocontrast}. Then the combination is regular semisimple, its centralizer in $\mathfrak{sp}(2n,\Qp)$ is a Cartan subalgebra of dimension $n$, and the $n$ independent Hessians, which commute with it, span this centralizer, so that definition implies the one used here. With only $n$ distinct eigenvalues the combination need not be regular, and the Cartan condition can fail.

For $n=2$, in coordinates with $\{z_i,w_j\}=\delta_{ij}$, the functions $q_1=z_1w_1+z_2w_2$ and $q_2=z_2w_1$ Poisson-commute and have independent Hessians, and the Hamiltonian matrix of $q_1+aq_2$ has exactly the $n=2$ distinct eigenvalues $\pm1$, but the Hamiltonian matrix of $q_2$ is nilpotent, so the span of $q_1,q_2$ is not a Cartan subalgebra. Crespo and Pelayo also compute $p$-adic Williamson normal forms effectively and show that they are independent of the real one \cite{crespopelayocontrast}; Theorem~\ref{thm:ev} turns such a normal form into a symplectic normalization of the system near the point.

The counterpart of the torus action comes from the Hamiltonians $\mathcal G^{-1}\circ \mathcal H$. They Poisson-commute and have linear flows $(z,w)\mapsto(e^{t\lambda}z,e^{-t\lambda}w)$ in the coordinates $\Phi^{-1}$, with weights $\lambda$. Take $l$ large enough for the exponentials to converge for $t\in p^l\Zp^n$. The resulting action of $p^l\Zp^n\simeq\Zp^n$ preserves $\mathcal H$.

\section{Logarithmic eigenvalues and height certificates}\label{sec:formal-eigen}

I now take linear parts whose eigenvalues are logarithms of algebraic points in the kernel of reduction of $\mathbb{G}_m$ or of an elliptic curve. Their small divisors are again logarithms of points, so the Smith exponents can be read off from exact arithmetic in the group, and heights bound them. Such eigenvalues lie between the algebraic case (Corollary~\ref{cor:smith-algebraic}) and the arbitrary case (Remark~\ref{rem:smith-liouville}). They arise for the generators of algebraic maps near rational fixed points (Remark~\ref{rem:log-origin}).

Let $\mathcal F$ be a one-dimensional commutative formal group law over $\Zp$ and $\omega(z)\,dz$ its normalized invariant differential, and let $\log_{\mathcal F}(z) = \int_0^z\omega$ be its logarithm \cite[Sections IV.4--IV.5]{silverman}. For $\mathbb{G}_m$, $\mathcal F = X + Y + XY$ and $\log_{\mathcal F}(z) = \log(1+z)$. For a minimal Weierstrass equation, the parameter $z = -x/y$ identifies the kernel of reduction with $p\Zp$ under the group law $\mathcal F$ \cite[Propositions VII.2.1--VII.2.2]{silverman}.

Let $p$ be odd, $k$ a number field and $\mathfrak p\mid p$ a prime of $k$ with $k_{\mathfrak p}=\Qp$, through which I embed $k$ in $\Qp$. Let $G$ be $\mathbb{G}_m$, or an elliptic curve given by a Weierstrass equation with coefficients in $k$ that is minimal at $\mathfrak p$. The local region is the kernel of reduction $G_1$: for $\mathbb{G}_m$ the $\alpha\in k$ with $\vp(\alpha-1)\ge1$, with parameter $z(\alpha)=\alpha-1$. For $E$ it consists of the points with $\vp(x)<0$, with $z=-x/y$ and $\vp(z)=-\vp(x)/2$, and $z(O)=0$. On $G_1$ put $\log_G=\log_{\mathcal F}\circ z$, that is, $\log_p\alpha$ or the $p$-adic elliptic logarithm $\log_E$. Since $p$ is odd, $\log_G$ is a homomorphism from $G_1$ to $p\Zp$ with $\vp(\log_GQ)=\vp(z(Q))$ \cite[Theorem IV.6.4 and Proposition VII.2.2]{silverman}. In particular it is injective and $G_1$ is torsion-free. I write the group law of $G$ additively.

\begin{proposition}[Logarithmic eigenvalues]\label{prop:log-eigen}
Let $P_1,\dots,P_d\in G(k)\cap G_1$, $\lambda_i=\log_G(P_i)$, and $A=S\Lambda S^{-1}$ with $\Lambda=\operatorname{diag}(\lambda_1,\dots,\lambda_d)$ and $S\in\mathrm{GL}_d(\Zp)$. For $\abs m=n$ and $1\le i\le d$ put $Q_{m,i}=\sum_jm_jP_j-P_i\in G(k)\cap G_1$. Then:
\begin{enumerate}[label=(\alph*)]
\item $\ip m\lambda-\lambda_i=\log_G(Q_{m,i})$, and the Smith exponents of $L_n$ are the $\vp(z(Q_{m,i}))$, with multiplicity;
\item an exponent is infinite exactly when $Q_{m,i}=O$ in $G(k)$; if $P_1,\dots,P_d$ are linearly independent in $G(k)$, then no $L_n$, $n\ge2$, has a zero Smith factor;
\item the exponents are $\vp(\alpha^m\alpha_i^{-1}-1)$ for $\mathbb{G}_m$ (with $P_j=\alpha_j$) and $\tfrac12\vp(1/x(Q_{m,i}))$ for $E$, computed exactly in $k$, without approximating any $\lambda_i$.
\end{enumerate}
\end{proposition}

\begin{proof}
(a) $\log_G$ is a homomorphism, and by Proposition~\ref{prop:smith-semisimple} with $S$ integral the Smith exponents of $L_n$ are the valuations of the eigenvalues $\ip m\lambda-\lambda_i$ of $L_\Lambda$, infinite for the zero ones. Since $\vp(\log_GQ)=\vp(z(Q))$ on $G_1$, they are the $\vp(z(Q_{m,i}))$. (b) $\log_G$ is injective on $G_1$. A relation $\sum_jm_jP_j=P_i$ with $\abs m=n\ge2$ is non-trivial, since $m\ne\varepsilon_i$. (c) restates the definitions.
\end{proof}

The $\lambda_i$ are in general not algebraic, so Corollary~\ref{cor:smith-algebraic} does not apply.

Heights now bound the exponents, as the Hadamard bound after Proposition~\ref{prop:smith-minors} does for rational matrices and the norm bound \eqref{eq:smith-algebraic} for algebraic eigenvalues. Let $h$ be the absolute logarithmic Weil height and $d_k=[k:\mathbb Q]$. On $E(k)$ let $\hat h(Q)=\lim_j4^{-j}h(x([2^j]Q))$ be the N\'eron--Tate height, a positive semi-definite quadratic form, and $\mu_E$ a constant with $h(x(Q))\le\hat h(Q)+\mu_E$ for all $Q\ne O$. This $\hat h$ is twice the canonical height of \cite[Section VIII.9]{silverman}, where its properties are proved. Theorem~1.1 of \cite{silverman90} bounds the difference between that canonical height and $\tfrac12h(x(Q))$ explicitly in terms of the heights of the $j$-invariant and of the discriminant, over any number field, and so gives an explicit $\mu_E$. Write $\Gamma$ for the Gram matrix of $\hat h$ on $P_1,\dots,P_d$, so that $\hat h(\sum_jm_j^\prime P_j)=(m^\prime)^{\mathsf T}\Gamma m^\prime$, and $\lambda_{\max}(\Gamma)$ for its largest eigenvalue.

\begin{proposition}[Height certificate]\label{prop:log-height}
In Proposition~\ref{prop:log-eigen}, every finite Smith exponent of $L_n$ satisfies
\begin{align}
s&\le\frac{d_k\bigl((n+1)h_*+\log2\bigr)}{\log p},\qquad h_*=\max_jh(\alpha_j),&&\text{for }\mathbb{G}_m,\label{eq:log-height-gm}\\
s&\le\frac{d_k\bigl(\lambda_{\max}(\Gamma)(n+1)^2+\mu_E\bigr)}{2\log p},&&\text{for }E.\label{eq:log-height-e}
\end{align}
Thus $\mu_n=O(n)$ for $\mathbb{G}_m$ and $\mu_n=O(n^2)$ for elliptic curves, with explicit constants and without any transcendence input.
\end{proposition}

\begin{proof}
If $\beta\in k$ has $\vp(\beta)=-e<0$, then $d_kh(\beta)\ge e\log p$, since $\mathfrak p$ has local degree one and contributes $\log\max(1,\abs\beta_{\mathfrak p})=e\log p$ to $d_kh(\beta)=\sum_vn_v\log\max(1,\abs\beta_v)$. Let $m^\prime=m-\varepsilon_i$, so that $\sum_j\abs{m_j^\prime}\le n+1$. For $\mathbb{G}_m$ put $\gamma=\prod_j\alpha_j^{m_j^\prime}\ne1$ and $\beta=1/(\gamma-1)$, so $e=s$. Then $h(\beta)=h(\gamma-1)\le h(\gamma)+\log2\le\sum_j\abs{m_j^\prime}h(\alpha_j)+\log2$, which gives \eqref{eq:log-height-gm}. For $E$ take $\beta=x(Q_{m,i})$, so $e=2s$, and $h(\beta)\le (m^\prime)^{\mathsf T}\Gamma m^\prime+\mu_E\le\lambda_{\max}(\Gamma)\norm{m^\prime}_2^2+\mu_E$ with $\norm{m^\prime}_2\le\sum_j\abs{m_j^\prime}\le n+1$.
\end{proof}

Yu's lower bounds for linear forms in $p$-adic logarithms improve \eqref{eq:log-height-gm} to a bound polynomial in $\log n$ \cite{yu}, with large constants depending on $p$, $k$ and the $\alpha_j$. Lower bounds for linear forms in $p$-adic elliptic logarithms are studied in \cite{remondurfels}. Proposition~\ref{prop:log-height} is elementary and explicit, and suffices for the following certificates.

\begin{corollary}[Certified exponents and a linear initial segment]\label{cor:log-certificate}
Suppose $P_1,\dots,P_d$ are linearly independent in $G(k)$, and let $B_n$ be the right-hand side of \eqref{eq:log-height-gm} or \eqref{eq:log-height-e}.
\begin{enumerate}[label=(\alph*)]
\item Every factor of the Smith decomposition of $L_n$ computed modulo $p^W$ with $W>B_n$ is non-zero, and the computed exponents are exact. In particular the directions kept as ``zero modulo $p^W$'' in Section~\ref{sec:normal} do not occur.
\item With threshold $\sigma$, $R_n=0$ in every degree $n$ with $B_n\le\sigma$. If $B_2\le\sigma$, this gives linearity up to degree $n_\sigma$, with
\[
n_\sigma+1=
\begin{cases}
\left\lfloor(\sigma\log p/d_k-\log2)/h_*\right\rfloor, & \mathbb G_m,\\
\left\lfloor\sqrt{(2\sigma\log p/d_k-\mu_E)/\lambda_{\max}(\Gamma)}\right\rfloor, & E.
\end{cases}
\]
If $B_2>\sigma$, the estimate gives no such degree $n\ge2$.
\end{enumerate}
\end{corollary}

\begin{proof}
By Proposition~\ref{prop:log-eigen}(b) all exponents are finite, and by Proposition~\ref{prop:log-height} they are at most $B_n<W$; Gaussian elimination modulo $p^W$ with pivots of least valuation then finds them exactly. For (b), all directions are removed when $B_n\le\sigma$. The bound $B_n$ increases with $n$.
\end{proof}

\begin{remark}[Where logarithmic eigenvalues come from]\label{rem:log-origin}
Let $p$ be odd, and let $g$ be a polynomial map with rational coefficients that is the time-one map $\varphi_1$ of a field $f$ satisfying \eqref{eq:hyp}. Time-one maps are congruent to the identity modulo $p^c$, coefficientwise, by Lemma~\ref{lem:coeff}. Conversely, if $g\equiv\mathrm{id}\pmod{p^c}$ coefficientwise with $c>1/(p-1)$, $g$ is the time-one map of the field $\sum_{j\ge1}(-1)^{j-1}\Delta_g^j(\mathrm{id})/j$, with $\Delta_g h=h\circ g-h$, which satisfies \eqref{eq:hyp}; its flow is Poonen's interpolation of the iterates of $g$ \cite[Proposition~3.8]{terskiiA}, \cite{poonen} (see also \cite{bell} and \cite[Theorem~2.1 and (A.14)]{cantatjang}). A fixed point $x^*$ of $g$ is an equilibrium of $f$, since $\norm{\varphi_1(x)-x}=\norm{f(x)}$ by Lemma~\ref{lem:orbit-metric}. The variational equation gives $Dg(x^*)=\exp(Df(x^*))$. If $x^*$ is rational and the eigenvalues $\alpha_i$ of $Dg(x^*)$ lie in $\Qp$, they are algebraic and lie in $1+p\Zp$, and the eigenvalues of $Df(x^*)$ are $\log_p\alpha_i$. So the normal forms of the generator of an algebraic map near a rational fixed point have $\mathbb{G}_m$-logarithmic eigenvalues. If moreover $Df(x^*)=S\Lambda S^{-1}$ with $S\in\mathrm{GL}_d(\Zp)$, which excludes Jordan blocks and non-integral eigenbases, then Propositions~\ref{prop:log-eigen} and~\ref{prop:log-height} apply with $h_*=\max_ih(\alpha_i)$, the heights of the multipliers of $g$; for a diagonalizable $Df(x^*)$ with $S\notin\mathrm{GL}_d(\Zp)$, Proposition~\ref{prop:smith-semisimple} adds the loss $(n+1)\kappa(S)$ to these bounds.
\end{remark}

\section{Zeros along an orbit of a linear normal form}\label{sec:formal-zeros}

To count zeros, I use linear normal coordinates and take the function to be polynomial in those coordinates. Let $L/\Qp$ be finite with ring of integers $\mathcal O_L$, $\vp(p)=1$, and suppose that near an equilibrium there are analytic coordinates $w$ in which the flow is $w(t)=e^{\Lambda t}w(0)$, with $\Lambda=\operatorname{diag}(\lambda_1,\dots,\lambda_d)$ and $\lambda_i\in p^c\mathcal O_L$ for an integer $c$ satisfying \eqref{eq:hyp}. Such coordinates exist in dimension one (the start of Section~\ref{sec:normal}), and in any dimension when the exact normal form of Section~\ref{sec:normal} has $R=0$, with $A=S\Lambda S^{-1}$, $S\in\mathrm{GL}_d(\mathcal O_L)$ and $w=S^{-1}\Phi^{-1}(z)$. The certificate of Theorem~\ref{thm:certificate} does not provide them, since it controls the normal form only modulo $p^K$. For an analytic $g=\sum_mg_mw^m$ and $w_0=w(0)$,
\[
G(t)=g(w(t))=\sum_\mu C_\mu e^{\mu t},\qquad C_\mu=\sum_{\ip m\lambda=\mu}g_mw_0^m,
\]
where $\mu$ runs over the distinct values of $\ip m\lambda$, all in $p^c\mathcal O_L$, and $C_\mu\to0$.

\begin{lemma}[Independence]\label{lem:independence}
Let the $\mu\in p^c\mathcal O_L$ be pairwise distinct and $C_\mu\in L$ with $C_\mu\to0$. Then $\sum_\mu C_\mu e^{\mu t}$ vanishes on $\Zp$ only if every $C_\mu$ is zero.
\end{lemma}

\begin{proof}
Since $\vp(C_\mu\mu^n/n!)\ge\vp(C_\mu)+n(c-\frac1{p-1})$ and $C_\mu\to0$, the terms $C_\mu\mu^nt^n/n!$ tend to zero jointly in $(\mu,n)$, uniformly in $t\in\Zp$, so the double series can be rearranged into the power series $\sum_nt^n\bigl(\sum_\mu C_\mu\mu^n\bigr)/n!$, which converges on $\Zp$. If it vanishes on $\Zp$, Strassmann's theorem \cite{strassmann} makes every coefficient zero, so $\sum_\mu C_\mu P(\mu)=0$ for every polynomial $P$. The functional $\phi\mapsto\sum_\mu C_\mu\phi(\mu)$ on $C(p^c\mathcal O_L,L)$ is continuous, and polynomials are dense there by Kaplansky's theorem \cite{kaplansky}, so it vanishes on the indicator of every ball $\B(\mu_0,\varepsilon)$. Given $\delta>0$, only finitely many $\mu$ have $\abs{C_\mu}\ge\delta$. Take $\varepsilon$ smaller than the distance from $\mu_0$ to those of them that differ from $\mu_0$. Then $0=\sum_{\mu\in\B(\mu_0,\varepsilon)}C_\mu$ is $C_{\mu_0}$ plus terms of norm less than $\delta$, so $\abs{C_{\mu_0}}<\delta$. As $\delta$ is arbitrary, $C_{\mu_0}=0$.
\end{proof}

Let $\Lambda_\lambda=\{k\in\mathbb Z^d:\ip k\lambda=0\}$ be the lattice of integer relations among the eigenvalues.

\begin{proposition}[Branches inside a subvariety]\label{prop:criterion}
Let $I$ be a set of analytic functions in the coordinates $w$. The whole orbit $\{w(t):t\in\Zp\}$ lies in the common zero set of $I$ if and only if $\sum_{m'\in m+\Lambda_\lambda}g_{m'}w_0^{m'}=0$ for every $g\in I$ and every class $m+\Lambda_\lambda$.
\end{proposition}

\begin{proof}
$\ip m\lambda=\ip{m'}\lambda$ exactly when $m-m'\in\Lambda_\lambda$, so these sums are the coefficients $C_\mu$ of $g(w(t))$, and Lemma~\ref{lem:independence} applies to each $g$.
\end{proof}

The kernels of $L_n$ detect only relations of the form $k=m-\varepsilon_i$, with $m\in\mathbb N^d$ and $\abs m=n$. For odd $p$, $\Lambda=\operatorname{diag}(2p,3p)$, $w_0=(1,1)$ and $g=w_1^3-w_2^2$ the orbit $(e^{2pt},e^{3pt})$ lies in $\{g=0\}$ because $(3,-2)\in\Lambda_\lambda$, although there is no resonance and every $L_n$, $n\ge2$, is invertible over $\Qp$.

\begin{proposition}[Zeros along an orbit]\label{prop:zeros}
Let $g$ be a polynomial in $w$, and let $\mu_1,\dots,\mu_s$, $s\ge2$, be the distinct values of $\ip m\lambda$ over the $m$ with $g_m\ne0$; if $g$ has degree at most $\delta$, then $s\le\binom{\delta+d}{d}$. Then for every $w_0$ with $G\not\equiv0$,
\[
\#\{t\in\Zp:G(t)=0\}\le s-1+\frac{s-2}{(p-1)c-1}\le 2s-3,
\]
and the number of zeros is at most $s-1$ when $(p-1)c>s-1$.
\end{proposition}

\begin{proof}
As in van der Poorten's argument \cite{vdpoorten} and \cite[Lemma~2]{vdpoortenschlickewei}, compare the Taylor coefficients by first looking at the derivatives. Write $G=\sum_na_nt^n$ and $b_n=n!a_n=\sum_jC_j\mu_j^n$, with $C_j=C_{\mu_j}$. The characteristic polynomial $\prod_j(x-\mu_j)=x^s+f_1x^{s-1}+\dots+f_s$ gives $b_n=-\sum_{k=1}^s f_kb_{n-k}$ for $n\ge s$. Its coefficients satisfy $\vp(f_k)\ge kc$ because $\mu_j\in p^c\mathcal O_L$. Hence $u_n=\vp(b_n)-cn$ satisfies $u_n\ge\min_{1\le k\le s}u_{n-k}$ and is bounded below by $m=\min_{k<s}u_k$ for all $n$. The minimum is finite because otherwise the recurrence gives $G\equiv0$.

Choose $k_0<s$ with $u_{k_0}=m$. Then $\vp(a_{k_0})\le m+c(s-1)$, while $\vp(a_n)\ge m+cn-(n-1)/(p-1)$ for $n\ge1$. If the Strassmann index $\nu$ is positive, comparing these two bounds gives $\nu(c-1/(p-1))\le c(s-1)-1/(p-1)$. Strassmann's theorem bounds the number of zeros by $\nu$ \cite{strassmann}, so this proves the first inequality, including the case $\nu=0$. The second uses $(p-1)c-1\ge1$ from \eqref{eq:hyp}.
\end{proof}

Given $s$, the bound depends only on $p$ and $c$. In particular it is uniform in $w_0$ and in the spacing of the exponents. For $p>s$ it is $s-1$, the bound for real exponential sums with $s$ terms. The proof uses only the linear recurrence satisfied by the derivatives of $G$ at $0$, so it applies unchanged to an orbit that is already written as a finite sum $\sum_\mu P_\mu(t)e^{\mu t}$ with polynomials $P_\mu$, with $s$ replaced by $\sum_\mu(\deg P_\mu+1)$. A resonant normal form can lead to a different kind of function. For odd $p$ the field $x'=px^2$ has linear part $0$ and a resonant quadratic term, and its orbit $x(t)=x_0/(1-px_0t)$ is not of this form. An algebraic subvariety in the original coordinates may become analytic in the linearizing coordinates. Using the finite-sum bound after truncation then requires a lower bound on the relevant coefficients $C_\mu$.

For integer times there is a bound of another kind: the number of integer zeros of a non-degenerate linear recurrence of order $s$ over a field of characteristic zero is bounded by a function of $s$ alone, independent of the coefficients and the roots \cite{schmidt,ess}. Proposition~\ref{prop:zeros} counts the zeros in $\Zp$ of the interpolating function, including the non-integer ones, which are the zeros that a $p$-adic certificate has to exclude. They occur already for $s=2$: since $\exp$ is injective on $p^c\mathcal O_L$, the function $C_1e^{\mu_1t}+C_2e^{\mu_2t}$ has exactly one zero in $\Zp$ when $-C_2/C_1$ lies in $\exp((\mu_1-\mu_2)\Zp)$ and none otherwise, and this zero is in general not an integer.

\section{Zero sets of orbits}\label{sec:zeros}

I now use modular and archimedean information to certify the full set of integer times at which an orbit of a polynomial map meets a hypersurface.

Let $F\colon\mathbb Z^d\to\mathbb Z^d$ be a polynomial map, $x_0\in\mathbb Z^d$, $x_n=F^n(x_0)$, and $g\in\mathbb Z[x_1,\dots,x_d]$. The question is the zero set $Z=\{n\ge0: g(x_n)=0\}$, or $\{n\in\mathbb Z: g(x_n)=0\}$ when the orbit runs backwards in $\mathbb Z^d$, as for a linear recurrence or a H\'enon map. For a linear map and a linear $g$ this is the Skolem problem, decidable up to order four (see \cite{bilu} and the references there) and, beyond that, for simple recurrences under the Skolem conjecture and the $p$-adic Schanuel conjecture \cite{bilu}. For nonlinear maps the dynamical Mordell--Lang conjecture predicts that $Z$ is a finite union of arithmetic progressions. The conjecture is known for \'etale maps and in several other cases \cite{bgt}, and no decision procedure is known. For linear recurrences, Bacik, Ouaknine, Purser and Worrell decide whether the $p$-adic interpolating function has a zero and compute the zeros exactly, with termination conditional on the $p$-adic Schanuel conjecture \cite{bacik}, and for maps $x+p^q\Phi(x)$ Mukhamedov gives finite certificates for the Strassmann index of the interpolating functions \cite{mukhamedov}. The sieve certifies $Z$ once its records rule out further zeros. Proposition~\ref{prop:nilpotent} gives a class of orbits on which this must happen.

A record collects what a single calculation can say about the zero set past a finite initial interval.

\begin{definition}\label{def:record}
A \emph{record} is a tail $T\ge0$ and a finite set $C$ of residue classes $c\bmod Q$, the cover, such that every $n\in Z$ with $n\ge T$ lies in a class of $C$, together with, for each class, a bound $b\in\N\cup\{\infty\}$ on the number of $n\in Z$ with $n\ge T$ in that class.
\end{definition}
A record with an empty cover says that $Z\subset[0,T)$. Two classes $c\bmod Q$ and $c'\bmod Q'$ meet if and only if $c\equiv c'\pmod{\gcd(Q,Q')}$, so covers are compared exactly. Records come from three sources.

A Strassmann record uses a prime where the cycle of the orbit modulo $p$ has an invertible Jacobian. In the chart of the cycle a power of the map is congruent to the identity coefficientwise (modulo $4$ for $p=2$), and the interpolation of Bell and Poonen \cite{bell,poonen}, applied to that power (see also \cite[Proposition~3.8]{terskiiA}), turns $g(x_n)$, on each class $n\equiv c\pmod Q$ past the tail, into a power series in $u$, with $n=c+Qu$. For a fixed $k$, the cover is the set of classes on which these series can vanish modulo $p^k$, and $b$ is Strassmann's bound for the class.

The orbit modulo $p^k$ gives a depth record just by inspecting the values of $g$ on its cycle.

\begin{lemma}\label{lem:depth}
Let $p$ be a prime and $k\ge1$, and let the orbit modulo $p^k$ have tail $T$ and period $P$. The classes $c\bmod P$ with $T\le c<T+P$ and $g(x_c)\equiv0\pmod{p^k}$, each with $b=\infty$, form a record. If $F$ induces a bijection of $(\mathbb Z/p^k)^d$ and the backward orbit lies in $\mathbb Z_{(p)}^d$, the orbit modulo $p^k$ is purely periodic and these classes contain every $n\in\mathbb Z$ with $g(x_n)=0$.
\end{lemma}

\begin{proof}
For $n\ge T$, $x_n\equiv x_c\pmod{p^k}$ with $c\equiv n\pmod P$ in $[T,T+P)$, and $g(x_n)=0$ gives $g(x_c)\equiv0$. A bijection of a finite set has only cycles, and then $x_n\equiv x_c$ for every $n\in\mathbb Z$ with $n\equiv c\pmod P$.
\end{proof}

Both hold, for instance, for the companion map of a linear recurrence whose last coefficient is prime to $p$, and for a map with a polynomial inverse over $\mathbb Z$, such as a H\'enon map.

An archimedean estimate can give a real record with an empty cover. A dominant root of a linear recurrence is one way to obtain such an estimate.

\begin{proposition}\label{prop:verdict}
Let $R_1,\dots,R_s$ be records with tails $T_j$, let $W\ge\max_j T_j$, and let $Z_W=Z\cap[0,W]$ be computed exactly. Suppose that for some $i$, every class of the cover of $R_i$ that meets a class of the cover of each $R_j$, $j\ne i$, has $b<\infty$ and contains at least $b$ elements $n\ge T_i$ of $Z_W$. Then $Z\cap[0,\infty)=Z_W$.
\end{proposition}

\begin{proof}
Let $n\in Z$ with $n>W$. Then $n\ge T_j$ for every $j$, so $n$ lies in a class of each cover, and these classes meet, since all contain $n$. The class of $R_i$ containing $n$ therefore has $b<\infty$ and already holds $b$ elements of $Z_W$ that are at least $T_i$. With $n$ these are $b+1$ elements of $Z\cap[T_i,\infty)$ in one class, which contradicts the bound.
\end{proof}

For an orbit that runs backwards the same argument is applied to $n\mapsto-n$: the records of Lemma~\ref{lem:depth} at primes where $F$ is bijective modulo $p^k$ are reflected, and the other records are computed for the backward orbit. A linear recurrence $u_n=a_0u_{n-1}+\dots+a_{k-1}u_{n-k}$ with $a_{k-1}\ne0$ runs backwards as the integer recurrence $s_m=a_{k-1}^m u_{-m}$, whose coefficients are $-a_{k-2-j}a_{k-1}^j$ for $j<k-1$ and $a_{k-1}^{k-1}$.

The adelic sieve by branch type (ABS) starts with the cheapest records and collects them until Proposition~\ref{prop:verdict} applies. At a prime $p$, compute the cycle of the orbit modulo $p$, its length $\ell$, and the product $J_{\mathrm{cyc}}$ of the Jacobians of $F$ along it, taken modulo $p$. The branch is \emph{\'etale} if $J_{\mathrm{cyc}}$ is invertible, \emph{nilpotent} if $J_{\mathrm{cyc}}^d\equiv0$, and \emph{mixed} if both parts of the Fitting decomposition are nonzero. Strassmann records require an \'etale branch. Depth records exist at every prime, for $k=1,2,\dots$. At \'etale primes they can exclude classes more cheaply. The counts $b$ come from Strassmann's theorem.

For a recurrence, the records at primes not dividing the last coefficient $a_{k-1}$ apply in both directions by Lemma~\ref{lem:depth}. At the other primes the two directions are treated separately. As records are added, the window $W$ grows to the largest tail. For a nilpotent branch, the orbit approaches a finite cycle in $\Zp^d$. The depth records eventually have empty cover precisely when $g$ is nonzero at every point of that cycle.

\begin{proposition}\label{prop:nilpotent}
Suppose that the orbit modulo $p$ has tail $T_1$ and a cycle of length $\ell$ along which the product $J_{\mathrm{cyc}}$ of the Jacobians of $F$ modulo $p$ is nilpotent. Then:
\begin{enumerate}[label=(\roman*),leftmargin=*]
\item For $0\le i<\ell$, the map $F^\ell$ has a unique fixed point $\xi_i$ in the residue disk of $x_{T_1+i}$, one has $\xi_{i+1}=F(\xi_i)$ (indices modulo $\ell$), and $x_{T_1+i+\ell m}\equiv\xi_i\pmod{p^{1+\lfloor m/d\rfloor}}$ for every $m\ge0$.
\item If $g(\xi_i)\ne0$ for every $i$, and $v_0=\max_i\vp(g(\xi_i))$, then the depth record with $k=v_0+1$ has an empty cover and a tail at most $T_1+\ell dv_0$. So every zero satisfies $n<T_1+\ell dv_0$.
\item If $g(\xi_i)=0$ for some $i$, the cover of every depth record at $p$ has a class that meets the class of $T_1+i$ modulo $\ell$.
\end{enumerate}
\end{proposition}

\begin{proof}
Let $D_i=x_{T_1+i}+p\Zp^d$. Since $F$ has integer coefficients, it maps $D_i$ into $D_{i+1}$, so $F^\ell$ maps $D_i$ into itself. The product of the Jacobians along the cycle starting at $x_{T_1+i}$ is a cyclic rotation of the product defining $J_{\mathrm{cyc}}$, hence nilpotent too, since $XY$ is nilpotent if and only if $YX$ is. Put $G=F^{\ell d}$. For $y,z\in D_i$, Taylor's formula with integral Hasse derivatives gives $G(z)-G(y)=DG(y)(z-y)+R$, where every term of $R$ has degree at least $2$ in $z-y$, so $\norm{R}\le p^{-1}\norm{z-y}$. By the chain rule $DG(y)$ is a product of $\ell d$ Jacobians at points congruent to the cycle modulo $p$, so it is congruent to the $d$-th power of a nilpotent $d\times d$ matrix, which is $0$ modulo $p$. Hence $\norm{G(z)-G(y)}\le p^{-1}\norm{z-y}$, and $G$ is a contraction of the complete space $D_i$ with a unique fixed point $\xi_i$. Since $F^\ell$ commutes with $G$, $F^\ell(\xi_i)\in D_i$ is fixed by $G$, so it equals $\xi_i$. Every fixed point of $F^\ell$ in $D_i$ is fixed by $G$, so $\xi_i$ is the only one, and $F(\xi_i)\in D_{i+1}$ is fixed by $F^\ell$, so it is $\xi_{i+1}$. From $x_{T_1+i}\equiv\xi_i\pmod p$ and the contraction, $\norm{x_{T_1+i+\ell dj}-\xi_i}\le p^{-1-j}$. Since $F^{\ell s}$ is $1$-Lipschitz and fixes $\xi_i$, this gives (i) for $m=dj+s$.

For (ii), take $n=T_1+i+\ell m$ with $m\ge dv_0$. By (i), $x_n\equiv\xi_i\pmod{p^{v_0+1}}$, so $g(x_n)\equiv g(\xi_i)\not\equiv0\pmod{p^{v_0+1}}$. So the orbit modulo $p^{v_0+1}$ is periodic from time $T_1+\ell dv_0$ on, and its tail is at most $T_1+\ell dv_0$. Every point of its cycle recurs at arbitrarily large times, so it is one of the $\xi_i\bmod p^{v_0+1}$, none of which is a zero of $g$ modulo $p^{v_0+1}$. The cover is empty.

For (iii), fix $k$. By (i), $x_n\equiv\xi_i\pmod{p^k}$ for every large $n\equiv T_1+i\pmod \ell$, so $g(x_n)\equiv0\pmod{p^k}$, and the class of such an $n$ past the tail is in the cover.
\end{proof}

ABS does not compute $\xi_i$ or $v_0$: it raises $k$ until the cover is empty, which by (ii) happens at the latest at $k=v_0+1$. The bound on the zeros is read off the computed cycle and involves no constants. Case (iii) needs a $p$-adic limit point of the orbit on the hypersurface $g=0$.

For $F(x,y)=(x^2,y^2)$, take $x_0=(2,0)$, $g(x,y)=x-4$ and $p=2$. The cycle modulo $2$ is the origin, with zero Jacobian and limit point $\xi=(0,0)$. Since $\vp(g(\xi))=2$, Proposition~\ref{prop:nilpotent} gives an empty cover at depth $3$. Modulo $8$ the orbit is $(2,0),(4,0),(0,0),\ldots$, so the tail has length two and the exact values $g(x_0)=-2$ and $g(x_1)=0$ give the zero set $\{1\}$. Along the orbit $y=0$, hence $\det DF=4xy$ vanishes there and no prime is \'etale.

To treat all three branch types at a prime $p$, consider the center manifold of the first return map. If $T_1$ and $\ell$ are the tail and cycle length modulo $p$, put $a=x_{T_1}$. In the coordinate $x=a+py$, the first return to its residue disk is
\[
\Psi(y)=\bigl(F^\ell(a+py)-a\bigr)/p.
\]
The coefficients of $\Psi$ are integers, and the coefficients in degree $e\ge1$ are divisible by $p^{e-1}$. Its linear part $\Lambda=D(F^\ell)(a)$ reduces to $J_{\mathrm{cyc}}$. Hensel's lemma applied to the characteristic polynomial of $\Lambda$ splits $\Zp^d=V_c\oplus V_s$, with both summands $\Lambda$-invariant, $\Lambda_c=\Lambda|_{V_c}$ invertible and $\Lambda_s=\Lambda|_{V_s}$ with $\Lambda_s^d\equiv0\pmod p$. The branch is \'etale for $V_s=0$ and nilpotent for $V_c=0$. In coordinates $(u,s)$ adapted to the splitting, with $u\in\Zp^\delta$, $\delta=\operatorname{rank}_{\Zp} V_c$, the linear part of $\Psi_c$ has no $s$-term and that of $\Psi_s$ no $u$-term. Gauss norms of series in $\Zp\langle u\rangle$ are suprema over $\mathcal O_{\mathbb C_p}^\delta$, so the estimates below may be checked pointwise there.

\begin{lemma}[Center manifold and shadow]\label{lem:shadow}
\begin{enumerate}[label=(\alph*),leftmargin=*]
\item There is a unique $\eta\in\Zp\langle u\rangle^{d-\delta}$ of Gauss norm at most $1$ whose graph $\mathcal A=\{s=\eta(u)\}$ satisfies $\Psi(\mathcal A)=\mathcal A$, and $\psi(u)=\Psi_c(u,\eta(u))$ is an isometric bijection of $\Zp^\delta$ with inverse in $\Zp\langle u\rangle^\delta$.
\item Every $y\in\Zp^d$ has a unique $y^*\in\mathcal A$ with $\norm{\Psi^m(y)-\Psi^m(y^*)}\le p^{-\lfloor m/d\rfloor}$ for all $m\ge0$.
\end{enumerate}
\end{lemma}

\begin{proof}
(a) For an integral series $\eta$ in the unit ball $\mathcal B$ of $\Zp\langle u\rangle^{d-\delta}$, the center component on its graph is $\psi_\eta(u)=\Psi_c(u,\eta(u))$. Its coefficientwise reduction is the invertible affine map $u\mapsto b_c+\Lambda_cu$. Solving for the inverse by contraction gives an isometric bijection on $\Zp^\delta$ and on $\mathcal O_{\mathbb C_p}^\delta$, with inverse in $\Zp\langle u\rangle^\delta$. Thus $\Psi$ sends the graph of $\eta$ to the graph of $\Gamma(\eta)=\Psi_s(\cdot,\eta(\cdot))\circ\psi_\eta^{-1}$, another element of $\mathcal B$.

Let $\epsilon=\norm{\eta_1-\eta_2}$. The dependence of $\Psi_c$ on $s$ is divisible by $p$, since its linear $s$-term vanishes and its terms of higher degree are divisible by $p$. It follows that the two maps $\psi_{\eta_i}$ and their inverses differ pointwise by at most $p^{-1}\epsilon$, using isometry for the inverses. The maps $\Psi_s$ and $\eta_2$ are $1$-Lipschitz, so the change of inverse contributes at most $p^{-1}\epsilon$ to the graph transform. The remaining leading difference is $\Lambda_s((\eta_1-\eta_2)\circ\psi_{\eta_1}^{-1})$, with error at most $p^{-1}\epsilon$. Compositions with these inverses preserve the Gauss norm. After $d$ steps the leading factor is $\Lambda_s^d\equiv0\pmod p$, so $\Gamma^d$ contracts by $p^{-1}$. Its unique fixed point $\eta$ is also fixed by $\Gamma$, since $\Gamma(\eta)$ is again fixed by $\Gamma^d$. Put $\psi=\psi_\eta$.

(b) Write $\tau(u,s)=s-\eta(u)$ for the difference from the graph and $\pi(u,s)=(u,\eta(u))$ for the projection onto it. Invariance of the graph gives
\[
\tau(\Psi(u,s))=\bigl[\Psi_s(u,s)-\Psi_s(u,\eta(u))\bigr]-\bigl[\eta(\Psi_c(u,s))-\eta(\Psi_c(u,\eta(u)))\bigr].
\]
Integral divided differences factor each bracket through $\tau$. For the first bracket the matrix is congruent to $\Lambda_s$ modulo $p$, since all higher degree terms are divisible by $p$. For the second the matrix is zero modulo $p$, by the vanishing of the linear $s$-term in $\Psi_c$. Hence $\tau(\Psi(z))=B(z)\tau(z)$ with $B\equiv\Lambda_s\pmod p$ and $\abs B\le1$. The distance from the graph cannot increase under one step, and every $d$ steps the product of these matrices reduces to $\Lambda_s^d=0$, giving a factor $p^{-1}$.

The center component also gives $\norm{\Psi(\pi z)-\pi\Psi(z)}\le\abs{\tau(z)}$, since distances between points on the graph are the distances between their $u$-coordinates. Pull $\pi\Psi^m(y)$ back along the isometry on the graph to $y^*_m$ with $\Psi^m(y^*_m)=\pi\Psi^m(y)$. Consecutive points satisfy $\norm{y^*_{m+1}-y^*_m}\le\abs{\tau(\Psi^m(y))}\le p^{-\lfloor m/d\rfloor}$, so their limit $y^*$ obeys $\norm{y^*_m-y^*}\le p^{-\lfloor m/d\rfloor}$. Taking the maximum of this distance and $\abs{\tau(\Psi^m(y))}$ gives the required bound. Two shadows would stay a constant distance apart under $\Psi$, so the distance bound forces them to coincide.
\end{proof}

\begin{proposition}\label{prop:branches}
For some positive multiple $M$ of $\ell$, each time class modulo $M$, represented by $r_0\in[T_1,T_1+M)$, has a series $G_{r_0}\in\Zp\langle t\rangle$ with
\[
\vp\bigl(g(x_n)-G_{r_0}(t)\bigr)\ge\lfloor(n-T_1)/(\ell d)\rfloor\qquad\text{for } n=r_0+Mt,\ t\ge0.
\]
\begin{enumerate}[label=(\roman*),leftmargin=*]
\item If $G_{r_0}$ has no zero in $\Zp$ and $\nu_{r_0}=\max_{t\in\Zp}\vp(G_{r_0}(t))$, then for every $k>\nu_{r_0}$ no class of the cover of the depth record at $p^k$ meets $r_0$.
\item If $G_{r_0}$ has a zero in $\Zp$, for every $k$ the cover of the depth record at $p^k$ has a class that meets $r_0$.
\end{enumerate}
The depth records at $p$ eventually have an empty cover exactly when every $G_{r_0}$ is zero-free on $\Zp$. In that case any $k>\max_{r_0}\nu_{r_0}$ suffices.
\end{proposition}

\begin{proof}
The reduction of $\psi$ is an invertible affine map, so a power $\phi=\psi^{M_1}$ is the identity modulo $p$ coefficientwise. For odd $p$ this power already meets the interpolation hypothesis, so put $M_2=M_1$. If $p=2$, write $\phi(u)=u+2h(u)$ with $h\in\Zp\langle u\rangle^\delta$. Substitution into this integral series gives $h(u+2h(u))\equiv h(u)\pmod 2$ coefficientwise, so
\[
\phi^2(u)=u+2h(u)+2h(u+2h(u))\equiv u\pmod 4.
\]
Taking $M_2=2M_1$ at $p=2$ means that in either case the coefficients of $\psi^{M_2}-\mathrm{id}$ have valuation greater than $1/(p-1)$. Poonen's theorem gives $\Phi\in\Zp\langle t,u\rangle^\delta$ with $\Phi(t,u)=\psi^{M_2t}(u)$ for $t\in\N$ \cite{poonen}. Put $M=\ell M_2$ and let $(w,\eta(w))$ be the shadow of $0$ from Lemma~\ref{lem:shadow}(b).

Each time class has a representative $r_0=T_1+i+\ell j$ with $0\le i<\ell$ and $0\le j<M_2$. Define $G_{r_0}(t)$ by evaluating $g\circ F^i$ at $a+p(U,\eta(U))$, with $U=\Phi(t,\psi^j(w))$ and the coordinates $(u,s)$ turned back into the original ones. The maps in this expression are integral, so $G_{r_0}$ belongs to $\Zp\langle t\rangle$. For $m=j+M_2t$, its value at an integer $t\ge0$ is $g(F^i(a+p\Psi^m(w,\eta(w))))$, while $x_{T_1+\ell m}=a+p\Psi^m(0)$. Applying the shadow estimate and the $1$-Lipschitz bounds for $g$ and $F^i$ gives $\vp(g(x_n)-G_{r_0}(t))\ge1+\lfloor m/d\rfloor\ge\lfloor(n-T_1)/(\ell d)\rfloor$ at $n=r_0+Mt$.

In (i), compactness and the absence of zeros make $\nu_{r_0}$ finite. Once $\lfloor(n-T_1)/(\ell d)\rfloor>\nu_{r_0}$, one has $\vp(g(x_n))=\vp(G_{r_0}(t))\le\nu_{r_0}$ throughout the class $r_0$. If a covered class $c\bmod P$ at depth $k>\nu_{r_0}$ met $r_0$, their intersection would contain arbitrarily large times by the Chinese remainder theorem. Past the modular tail such times would have $g(x_n)\equiv0\pmod{p^k}$, which contradicts the bound.

For (ii), fix $k$ and let $G_{r_0}(t_*)=0$ with $t_*\in\Zp$. By continuity there is an integer $e\ge0$ such that $\vp(G_{r_0}(t))\ge k$ throughout $t_*+p^e\Zp$. Choose an integer $t$ in it large enough that $n=r_0+Mt$ lies past the modular tail and satisfies $\lfloor(n-T_1)/(\ell d)\rfloor\ge k$. Then $g(x_n)\equiv0\pmod{p^k}$, so the cover contains a class meeting $r_0$.
\end{proof}

On a nilpotent branch, $\delta=0$ and each $G_{r_0}$ is the constant $g(\xi_i)$ for the corresponding point of the limit cycle. Proposition~\ref{prop:nilpotent} gives the depth bound in this case. On an \'etale branch, the center manifold is the whole chart and $G_{r_0}$ is the Bell--Poonen interpolation used in a Strassmann record. A $p$-adic zero of $G_{r_0}$ keeps the class in the cover at every depth. On an \'etale branch one can still close the class if the exact zeros already found exhaust its Strassmann bound. A non-integer zero prevents that count from being exhausted, so the sieve then needs a record from another place.

To make a depth record I follow the finite orbit modulo $p^k$. When a class survives, I can raise $k$ and repeat the calculation. The limit cycle and the shadow tell me when that will help, using the same kind of attraction that appears in Ghioca and Tucker's work on dynamical Mordell--Lang \cite{ghiocatucker}.

How often can the sieve make a Strassmann record? Benedetto, Ghioca, Hutz, Kurlberg, Scanlon and Tucker model the orbit of a point modulo $p$ by a random map of a set of $p^d$ points, with its cycle independent of the hypersurface $\det dF\equiv0$ \cite[Section~5]{bghkst}. A cycle of length $\ell$ then avoids the hypersurface, that is, the branch at $p$ is \'etale, with probability $(1-1/p)^\ell$. The cycle length is typically of size $p^{d/2}$, but short cycles dominate the average, and the probability $\pi_p$ that the branch is \'etale is $(\sqrt{\pi/2}+o(1))\,p^{1-d/2}\approx1.25\,p^{1-d/2}$ for $d\ge3$, and about $0.6$ for $d=2$ \cite[Proposition~5.3 and Remark~5.4]{bghkst}. Assuming independence across primes, they note that the probability of no \'etale prime is an Euler product, which vanishes for $d=3,4$ and is positive for $d\ge5$, and they test the model in dimensions up to three \cite[Sections~5.1.3 and~5.2]{bghkst}. The same model gives rates for the number of such primes and for the probability of finding none below a cutoff.

For a birational map that contracts a curve to a point, take a prime $p$ at which the reduction still contracts the curve to a point, is injective outside it, and has the curve as the full preimage of its image. If the starting point is not periodic modulo $p$, its eventual cycle modulo $p$ passes through the curve, where $\det dF$ vanishes modulo $p$, so the branch at $p$ is not \'etale, at every such prime at once. Birational maps without a contracted curve can do the opposite: the H\'enon map $(x,y)\mapsto(y,x+y^2)$ has a polynomial inverse and $\det dF=-1$, so its branch is \'etale at every prime. Either way such maps fall outside the random-map model, so the heuristic is for nonbirational maps. A point with an infinite orbit can still be periodic modulo every odd $p$: for $(x,y)\mapsto(x,xy)$ the point $(2,1)$ has a cycle of length $\ell=\operatorname{ord}_p2$ with multiplier $2^\ell$, so its branch is \'etale at every odd $p$.

\begin{heuristic}[after \cite{bghkst}]\label{heur:etale}
For random quadratic maps of $\mathbb Z^d$, with coefficients uniform in $[-H,H]$ and $H\to\infty$, that are not birational, and a random starting point with an infinite orbit, the events ``the branch at $p$ is \'etale'' are modeled as independent across primes, with probabilities $0\le\pi_p<1$ and $\pi_p\sim c_dp^{1-d/2}$, where $c_d=\sqrt{\pi/2}$ and $d\ge3$. Put $s=d/2-1$. The model gives the following predictions.
\begin{enumerate}[label=(\alph*),leftmargin=*]
\item the expected number of \'etale primes $p\le X$ grows like $2c_3\sqrt X/\log X$ for $d=3$ and like $c_4\log\log X$ for $d=4$, and stays bounded for $d\ge5$;
\item the logarithm of the probability that no prime $p\le X$ is \'etale is asymptotic to $-2c_3\sqrt X/\log X$ for $d=3$, and the probability itself is $(\log X)^{-c_4+o(1)}$ for $d=4$;
\item for $d\ge5$ the probability that no prime at all is \'etale is positive and equals $\prod_p(1-\pi_p)$, the limit as $P\to\infty$ of
\[
\prod_{p\le P}(1-\pi_p)\cdot\prod_{p>P}\bigl(1-c_dp^{-s}\bigr),\qquad \prod_{p>P}\bigl(1-c_dp^{-s}\bigr)=\Bigl(\zeta(s)\prod_{p\le P}(1-p^{-s})\Bigr)^{-c_d}\bigl(1+O(P^{1-2s})\bigr).
\]
\end{enumerate}
\end{heuristic}

The three regimes are those of $\sum_pp^{-s}$: below the pole of $\zeta$ at $s=1$ for $d=3$, at the pole for $d=4$, and past it for $d\ge5$. The dichotomy between $d\le4$ and $d\ge5$ is that of \cite{bghkst}. Partial summation gives the rate for $d=3$, and Mertens' theorem for primes gives the rate for $d=4$. For $d\ge5$ the sum of the probabilities converges, and expanding the logarithm of the tail product gives (c). In dimension four the sum diverges like $\log\log X$, so Bell's argument at one prime is expected to apply to almost every orbit, but the share of orbits without an \'etale prime up to $X$ decays only like $(\log X)^{-c_4}$.

So, under the independence in Heuristic~\ref{heur:etale}, Strassmann records, which need an \'etale prime, would be missing for a positive proportion of orbits from dimension five on, while depth records exist at every prime.

\section{Further questions}\label{sec:limitations}

At a point $m$ of rank $\ell>0$, let $V$ be the span of the Hamiltonian vectors at $m$, so that the transverse system lives on the symplectic quotient $V^\omega/V$. The classical Carath\'eodory--Jacobi--Lie induction uses local flows and the inverse function theorem and suggests taking $H_i-H_i(m)=\xi'_i$ for $i\le\ell$ as part of a system of symplectic coordinates. The remaining integrals would then give a family of rank-$0$ problems parametrized by $\xi'$, with the same integer divisors \eqref{eq:ev-divisor} for each parameter. A proof with a radius uniform in $\xi'$ would require analytic dependence of the critical point and eigencoordinates, as well as a Newton step in Gauss norms that include the parameter variables.

For the radius in Theorem~\ref{thm:ev}, I would start with the four-dimensional models of \cite{crespopelayoww}. Computing $r$ and $\mu$ there would show how the contribution from the eigencoordinates and $\Lambda$ changes with Williamson type and ramification. How sharp is the dependence on $p$ through $L_p$?

For a single field with semisimple non-resonant linear part, let $\omega(D)$ be the largest valuation of a divisor $\ip m\lambda-\lambda_i$ with $\abs m\le D$. In the integrable case the sum $\sum_k s_{N_k}/2^k$ controls the losses, which suggests $\sum_k\omega(2^{k+1})/2^k<\infty$ as a condition for linearization of a single field. The Liouville example of Remark~\ref{rem:smith-liouville} fails this Brjuno condition. Does the condition suffice, and is it sharp? Resonances pose a different problem for a Newton scheme: the bracket of a retained quadratic term with a degree-$N$ correction appears already in degree $N+1$. A quadratic version of the fixed-threshold construction would have to control those terms.

If $g(x_n)=0$ at $n=r_0+Mt$, Proposition~\ref{prop:branches} forces $\vp(G_{r_0}(t))\ge\lfloor(n-T_1)/(\ell d)\rfloor$. An infinite orbit can have $G_{r_0}\equiv0$, as happens for $F(x,y)=(3x,3y)$, $x_0=(1,0)$, $g=y$ and $p=3$. For nonzero $G_{r_0}$, Weierstrass preparation writes it as $p^\nu$ times a unit and a polynomial. If there are infinitely many integer zeros, some root $\tau$ must satisfy $\vp(t-\tau)\ge ct-C$ along a subsequence for fixed $c>0$ and $C$. The approximation puts $\tau$ in $\Zp$, and it cannot be an integer because the valuation of a nonzero integer difference grows at most logarithmically in $t$. Can algebraic $F$, $x_0$ and $g$ give such an approximation? This is akin to the transcendence questions used for simple recurrences in \cite{bilu}.

On a nilpotent branch whose limit cycle meets $g=0$, Proposition~\ref{prop:nilpotent}(iii) leaves classes in every depth record at that prime. Xie's theorem gives a finite union of arithmetic progressions for the full zero set in dimension two over $\bar{\mathbb Q}$ \cite{xie}, and the corresponding case in higher dimension remains open. Heuristic~\ref{heur:etale} raises a different question for a fixed orbit: how often does its reduction have an \'etale cycle? In dimension one this leads to Chebotarev questions for preimage towers and arboreal Galois representations. For a Zariski-dense orbit in dimension two, one may ask whether its \'etale primes have positive density.

\section{Code availability}\label{sec:code}

The \texttt{galrax.padic} package in Galrax (\url{https://github.com/NotTruecCoder/galrax}) implements the Smith normal forms and their certificates in \texttt{\fontencoding{T1}\selectfont normal\_form.py}. The Eliasson--Vey computation is \texttt{\fontencoding{T1}\selectfont IntegrableSystem.normal\_form} in \texttt{integrable.py}. The code for logarithmic eigenvalues is in \texttt{formal.py}. The records and the sieve are in \texttt{adelic.py}.

\end{document}